\documentclass[reqno]{amsart}
\usepackage{graphicx,amsmath, amssymb,amscd}
\renewcommand{\Im}{{\mathop{\mathrm{Im}}\nolimits}}
\renewcommand{\Re}{{\mathop{\mathrm{Re}}\nolimits}}
\renewcommand{\kappa}{\varkappa}
\renewcommand{\phi}{\varphi}

\newcommand{\eps}{\varepsilon}
\newcommand{\R}{{\mathbb R}}
\newcommand{\C}{{\mathbb C}}
\newcommand{\B}{{\mathbb B}}
\newcommand{\dg}{\mathop{\mathrm{diag}}\nolimits}
\newcommand{\Sym}{\mathop{\mathrm{Sym}}\nolimits}

\newcommand{\abs}[1]{\left|{#1}\right|}
\newcommand{\tps}[1]{{#1}^{t}}
\newcounter{romct}
\newenvironment{romanlist}{\begin{list}{(\roman{romct})}{\usecounter{romct}}}{\end{list}}

\newtheorem{theorem}{Theorem}[section]
\newtheorem{lemma}{Lemma}[section]
\newtheorem{proposition}{Proposition}[section]
\newtheorem{corollary}{Corollary}[section]
\newtheorem{definition}{Definition}[section]

\begin{document}
\title{Holomorphic extension of CR functions from quadratic cones}
\author{Debraj Chakrabarti\and Rasul Shafikov}
%
%

\maketitle
\begin{abstract}
It is proved that CR functions on a quadratic cone $M$ in $\C^n$,
$n>1$, admit one-sided holomorphic extension whenever $M$ does 
not have two-sided support, a geometric condition on $M$ which
generalizes minimality in the sense of Tumanov. A biholomorphic
classification of quadratic cones in $\C^2$ is also given.
\end{abstract}

\section{Introduction.} \label{intro}One of the central results in the
theory of CR functions is Tr\'epreau's theorem
\cite{trepreau:theorem} on local holomorphic extension of CR
functions defined on a smooth real hypersurface $M$ in $\C^n$ to one
side of $M$. Its generalization \cite{tumanov}, \cite{br} to the
case when $M$ is of higher codimension is known as
wedge-extendability. In both cases simultaneous extension of all CR
functions to an open set in $\C^n$ adjacent to $M$ is equivalent to
{\em minimality} of $M$. In the hypersurface case this simply means
that $M$ does not contain any germ of complex analytic hypersurface.
For more details on this subject see the recent extensive survey
\cite{mp} and the references there.

CR functions can be defined on a wider class of objects than smooth
CR submanifolds, for example on locally Lipschitz graphs (see e.g.,
\cite{chirka-stout}). It is therefore natural to study the
properties of CR functions defined on non-smooth objects, in
particular, a natural question to ask is whether a similar one-sided
holomorphic extension of CR functions holds when the hypersurface
$M$ is no longer smooth. More generally, one can ask whether all
functions holomorphic on one side of $M$ extend holomorphically to
the other side near non-smooth points of $M$. The answer is
affirmative if $M$ is the graph of a continuous function, as proved
in \cite{chirka:graphs}. While the problem could be very difficult
in general, it seems that the situation when $M$ is a real analytic
variety is of particular interest because of the connection with
questions related to boundary regularity and analytic continuation
of holomorphic mappings between domains with real analytic boundary.

In this paper we consider the family of real quadratic cones in
$\C^n$, $n>1$, and prove a Tr\'epreau-type extension result for CR
functions. By a {\em (real) quadratic cone} we mean an irreducible
real algebraic set $M$ in $\C^n$ of pure dimension $2n-1$ defined by
\begin{equation}\label{cone-def}
M=\{z\in\C^n\colon \rho(z)=0\},
\end{equation}
where $\rho(z)$ is a real-valued homogeneous polynomial  of degree
two in  $x_j$ and $y_j$, where $z_j=x_j+iy_j$, $1\leq j\leq n$, are
the coordinates on  $\C^n$. These cones are perhaps the simplest
examples of real analytic varieties with singularities, and they
form a natural class of models to study. Let $M^{\rm reg}$ be the
smooth part of $M$. We call $f\in L^1_{\rm loc}(M)$ a CR function,
if $\int_{M^{\rm reg}} f\overline\partial\phi=0$ for every smooth
$(n,n-2)-$form $\phi$ with compact support in $\C^n$. In particular,
every continuous function on $M$, which is CR on $M^{\rm reg}$, is
CR on $M$ (see Lemma~\ref{crexamples} below.) We refer to
Sections~2~and~3 for further discussion concerning quadratic cones
and CR functions. It turns out that minimality is no longer a
sufficient condition for one-sided extension. Indeed, consider
\begin{equation}\label{example-M}
M=\left\{z=(z_1,z_2)\in\C^2 : \rho(z) =
\Re\left(\frac{1}{2}z_1^2+\frac{1}{3}z_2^2\right) + |z_1|^2 -
|z_2|^2=0\right\}.
\end{equation}
Then $M$ is a quadratic cone of dimension three, in fact, $M$ is a
smooth minimal manifold away from the origin, and $\C^2\setminus M$
consists of two connected components $\Omega^\pm=\{z\in\C^2:
\pm\rho(z)>0\}$. Consider the function
\begin{equation}\label{example-f}
f(z)= \frac{z_2^2}{z_1}-\frac{z_1^2}{z_2}.
\end{equation}
By letting $f(0)=0$, it is easy to see that $f|_M$, the restriction
of $f$ to $M$,  extends continuously to the origin, and therefore
$f|_M$ is a continuous CR function. Like CR functions on smooth
hypersurfaces, $f$ admits the  jump representation $f=F^+ - F^-$,
where $F^+=\frac{z_2^2}{z_1}\in\mathcal O\left(\Omega^+\right)$ and
$F^-=\frac{z_1^2}{z_2}\in\mathcal O\left(\Omega^-\right)$. However,
since $\{z\in \C^2: z_2=0\}$ is contained in $\overline{\Omega^+}$,
and $\{z_1=0\}$ is contained in $\overline{\Omega^-}$, it follows
that $f|_{M}$ does not extend holomorphically to either side of $M$.

In the above example the obstruction to holomorphic extension is the
presence of two complex hypersurfaces $A^+=\{z_2=0\}$ and
$A^-=\{z_1=0\}$ residing on different sides of $M$. We prove that
for quadratic cones this is the only obstruction. To formulate this
condition we make the following definition. Given a real analytic
hypersurface $M$, i.e., a real analytic set of pure dimension $2n-1$
in $\C^n$, we call $\rho(z)$ a {\it defining function} of $M$ near
point $p\in M$ if $\rho$ is real analytic in a neighbourhood
$\Omega$ of $p$, $M\cap \Omega=\{z\in \Omega: \rho(z)=0\}$, and
there exists a point $z\in M\cap \Omega$ such that $d \rho(z)\ne 0$.
Clearly, any real analytic hypersurface $M$ admits a defining
function at each of its points.

Let $\Omega^\pm=\{z\in \Omega : \pm\rho(z)>0\}$. Then $\Omega^\pm$
are non-empty open sets which up to a sign are independent of the
choice of the defining function of $M$, as proved in
Lemma~\ref{omegas} below. We note that in general $\Omega^\pm$ may
consist of several connected components.

\begin{definition}\label{support}
  Let $M$ be an irreducible real analytic hypersurface in $\C^n$. Let
  $\rho(z)$ be a defining function of $M$ in some neighbourhood $\Omega$ of
  a point $p\in M$. We say that $M$ admits two-sided support by
  complex hypersurfaces at $p\in M$ if there exist germs at $p$ of
  complex analytic hypersurfaces $A^+$ and $A^-$ such that
  $A^\pm\subset \overline{\Omega^\pm}$.
\end{definition}

In short we say that $M$ as above has {\em two-sided support at
$p$}. As in the smooth case, we say that a real analytic
hypersurface $M$ is minimal at $p\in M$, if there is no germ of
complex analytic hypersurface at $p$ contained in $M$. Note, in
particular, that any hypersurface $M$, which is non-minimal at $p\in
M$, has two-sided support at $p$.

It is well-known (see \cite{kohnnirenberg:counterexample}) that a
smooth pseudoconvex hypersurface need not have a supporting complex
hypersurface on the pseudoconcave side, and therefore in general the
existence of such support is irrelevant to holomorphic
non-extendability of domains in $\C^n$. However, in the question of
simultaneous analytic continuation of all holomorphic functions from
one side of $M$ to the other side, the existence of supporting
complex hypersurfaces becomes important. Our first result concerns
the extension of holomorphic functions defined on one side of a
cone.

\begin{theorem}\label{extensionthm}
  Let $M=\{\rho=0\}$ be a quadratic cone in $\C^n$, $n>1$, where
  $\rho$ is the defining function of~$M$. Given a neighbourhood $\Omega$
  of $0$ in $\C^n$, let $\Omega^+$ and $\Omega^-$ be the open sets in $\C^n$
  defined by $\Omega^\pm=\{z\in\Omega\colon\pm\rho(z)>0\}$. Then the following
  are equivalent.

  (a) $M$ does not have two-sided support at the origin.

  (b) There exists a neighbourhood $U$ of the origin such that either
  all functions in $\mathcal O(\Omega^+)$ or all functions in $\mathcal
  O(\Omega^-)$ extend holomorphically to $U$.
\end{theorem}

We note that in the above theorem the singular part of the cone $M$
can have arbitrary dimension. Quadratic cones provide some
understanding of singularities of arbitrary analytic sets. We
believe that Definition \ref{support} is also relevant to extension
results for arbitrary real analytic sets. This is the subject of our
further investigation.

Combining Theorem~\ref{extensionthm} with the jump formula for CR
functions, we obtain the following result.

\begin{corollary}\label{cr-ext}
  If $M$ does not have two-sided support at the origin, then there
  exists a neighbourhood $U$ of the origin such that every CR function
  on $\Omega\cap M$ has a holomorphic extension either to $U\cap \Omega^+$
  or $U\cap \Omega^-$. In particular, if $f$ is a continuous function
  on $M\cap\Omega$ and CR on $M^{\rm reg}$, then there is a
  function~$F$, holomorphic in $U\cap \Omega^+$ or $U\cap \Omega^-$ and
  continuous on $M^{\rm reg}\cap U$, with $F|_{M^{\rm reg} \cap U}=f$.
\end{corollary}

It remains an open question whether for continuous $f$ the equality
$F=f$ holds on the singular part of $M\cap U$. An answer will depend
upon a better understanding of boundary behavior of the jump
formula, for example as in \cite{rosay}.

We also note here that in  the non-smooth case the well-known
approximation theorem due to Baouendi and Tr\`{e}ves \cite{bt},
fails in general. The conditions under which a CR function on a
non-smooth hypersurface can be approximated by holomorphic
polynomials are not known. In fact, for $M$ given by
\eqref{example-M} in the example above, not all CR functions on $M$
can be approximated by polynomials. Indeed, if this were true, then
by the maximum principle, approximation by polynomials would also
hold on the union of all holomorphic discs attached to $M$ near the
origin. Fix some small $\epsilon>0$ and let $D_\epsilon = A^+\cap
\{|z_1|\leq\epsilon\}$. Then $\partial D_\epsilon$ is compactly
contained in $\Omega^+$. If $\tilde{D}_\epsilon$ is a small
perturbation of $D_\epsilon$, the set $\tilde{D}_\epsilon\cap
\overline{\Omega^-}$ is a disc attached to $M$ whenever its interior
is nonempty. We can therefore produce a family of holomorphic discs
attached to $M$ that sweep a small one-sided neighbourhood of the
origin. But this would imply that any CR function can be extended
holomorphically to one side of $M$. However, the function in
\eqref{example-f} does not admit such extension.

Tr\'epreau's original proof \cite{trepreau:theorem} relies on the
existence of Bishop disks attached to a smooth hypersurface.
Alternatively one can use the technique of propagation of
holomorphic extendability along CR orbits (see \cite{mp}), or
Shcherbina's ~\cite{shcherbina} characterization of polynomial hulls
of continuous 2-spheres in $\C^2$ (as in \cite{chirka:graphs}). The
last technique yields the analog of Tr\'epreau's theorem even for
graphs of continuous functions. However, these methods require the
hypersurface either to be smooth or locally represented as a graph,
which is impossible for quadratic cones near the singular points.

We therefore follow a different route to prove
Theorem~\ref{extensionthm}.  For $n=2$, when $M$ does not have
two-sided support, we explicitly construct a family of analytic
discs touching $M$ at the origin and then apply the
Kontinuit\"{a}tssatz.  A similar technique was also used in
\cite{bf:mmj1978}. The case $n>2$ is then reduced to $n=2$ by using
biholomorphic invariants of cones, which will be obtained in the
proof of Theorem~\ref{normalform} below, and a slicing argument.

In order to construct the discs we obtain a complete classification
of quadratic cones in $\C^2$ under biholomorphic equivalence. Such
biholomorphic classification of non-smooth hypersurfaces is of
independent interest (see e.g., \cite{gongburns:leviflat} for
classification of non-smooth Levi-flats.)

\begin{theorem}\label{normalform}
Let $M$ be a quadratic cone. Then there is a unique normal form
$\rho$ in the last column of the table below, such that $M$ is
biholomorphic to $\{\rho=0\}$. More precisely, there is a unique
type in column~2 of the table determined by the signature of the
hermitian part of $\rho$ given in column~1, and a unique choice of
parameters in the range given in column~3 so that the germ of $M$ at
the origin is biholomorphic to the corresponding germ of
$\{\rho=0\}$, where $\rho$ is the defining function in the normal
form in the last column of this row with these parameters.
\end{theorem}

\begin{center}
\begin{tabular}{|p{1cm}|p{1.5cm}|p{2cm}|p{5cm}|}
  \hline
  $(\pi,\nu)$ & Type & Parameters & Defining Function \\\hline
  (2,0) & $\mathcal{M}_{(2,0)}$ & $\begin{array}{c}0\leq B\leq A,\\A>1\end{array}$
   & $\Re (A z_1^2 + B z_2^2) +
  |z_1|^2 +  |z_2|^2$ \\\hline
 & $\mathcal{M}_{(1,1)}^1$ & $0\le B \le
  A$ & $\Re (A z_1^2 + B z_2^2) +
  |z_1|^2 -  |z_2|^2$ \\\cline{2-4}
  (1,1) & $\mathcal{M}_{(1,1)}^2$ &$\begin{array}{c}
  \Re A > 0,\\ \Im A\ge 0\end{array}$  & $\Re(A z_1^2+ \overline{A} z_2^2)+\Im
  (z_1\overline z_2)$ \\\cline{2-4}
 & $\mathcal{M}_{(1,1)}^3$ &    & $\Re( z_1^2)+\Im(z_1\overline z_2)$ \\\cline{2-4}
  & $\mathcal{M}_{(1,1)}^4$ & $A>0$ & $\Re ( z_1^2 + iA z_1z_2) +
  \Im(z_1\overline{z}_2)$ \\
 \hline
 (1,0)   & $\mathcal{M}_{(1,0)}^1$ &$A\geq 0$ &$\Re(Az_1^2+z_2^2)+\abs{z_1}^2$ \\\cline{2-4}
   & $ \mathcal{M}_{(1,0)}^2$&  & $\Re(z_1z_2)+\abs{z_1}^2$ \\\hline
 (0,0) & $\mathcal{M}_{(0,0)}^1$  &  & $\Re(z_1^2+z_2^2)$ \\\cline{2-4}
  \hline
\end{tabular}
\end{center}
\bigskip
{\em Note:} After submitting this article for publication, we came
to know of the work of Ermolaev \cite{ermo}, in which real quadratic
forms in $\C^n$ are classified by an algebraic method very different
from ours. However, our method yields directly the invariants
required for the proof of Theorem~\ref{extensionthm} for $n\geq 3$.

Two-sided support by complex hypersurfaces is a natural
generalization of the notion of non-minimality. The following result
is probably well-known, but for the sake of completeness we will
give the proof.

\begin{proposition}\label{smoothconda} Let $M \subset\C^n$ be a $\mathcal{C}^1$
hypersurface, or more generally, suppose that $M$ can be represented
locally as a Lipschitz graph. Let $p\in M$, and let $M$ divide a
neighbourhood $\Omega$ of $p$ in $\C^n$ into two components
$\Omega^\pm$. Suppose that $M$ has two-sided support at $p$, so that
there are germs at $p$ of complex hypersurfaces $A^{\pm}\subset
\overline{\Omega^\pm}$. Then $A^+=A^-\subset M$.
\end{proposition}

The paper is organized as follows: in Section~2 we discuss
properties of quadratic cones, in Section~3 we review the material
on CR functions and prove Corollary~\ref{cr-ext}. Section~4 contains
some preliminary results in Linear Algebra. Theorem~\ref{normalform}
is proved in Section~5 and Theorem~\ref{extensionthm} in Section~6.
As an application of the proof of Theorem~\ref{extensionthm} we
classify in Section~7 all quadratic cones that have two-sided
support. Finally, Proposition~\ref{smoothconda} is proved in
Section~8.

{\em Note on this version:} An earlier version of this article was published as
\cite{mathann}, and subsequently an erratum \cite{erratum} correcting some mistakes was
published. This version incorporates the corrections from \cite{erratum}, as well as
corrections of a few minor typos, and is posted
in ArXiv.

{\em Acknowledgment.} The authors would like to thank A.~Boivin,
X.~Gong,\linebreak[4] J.-P.~Rosay and A.~Sukhov for helpful
discussions.

\section{Quadratic Cones.}

\subsection{Real analytic hypersurfaces have two sides.}
We say that $M$ is a real analytic subset of a domain $U\subset
\R^N$  if near each point $p\in U$ it is defined as the zero set of
a finite collection of real analytic functions near $p$. $M$ is
called irreducible if it cannot be represented as a union of two
non-empty real analytic sets. A point $p\in M$ is called {\it
regular} or smooth of dimension $d$, if near $p$, $M$ is a
$d$-dimensional real analytic manifold. Dimension of $M$ is then
defined as the maximum dimension at regular points. We say that $M$
has pure dimension $d$ if all regular points of $M$ have dimension
$d$. We will say that $M$ is a {\em real analytic hypersurface} in
$\R^N$ if it is of pure dimension $N-1$. Denote by $M^{\rm reg}$ the
set of regular points of $M$. Then $M^{\rm sng}=M\setminus M^{\rm
reg}$ is the set of singular points.

\begin{lemma}\label{omegas}
Let $M\subset \R^N$ be an irreducible real analytic hypersurface
which is given in some open set $\Omega$ by the defining function
$\rho$ (recall from Section~\ref{intro} that this means that there
is a point $x\in M\cap\Omega$ such that $d\rho(x)\ne 0$). Let
$\Omega^\pm=\{\pm\rho>0\}$. Then up to a sign, the sets $\Omega^\pm$
are independent of the choice of the defining function of $M$.
\end{lemma}

\begin{proof} Let $S=\{x\in\Omega\colon d\rho(x)=0\}\cap M$. Then $S$
  is a real analytic set of dimension at most $N-2$, in particular,
  $S$ is nowhere dense in $M$, and $M^{\rm sng}\subset S$. Suppose
  $\tilde{\rho}(x)$ is another defining function of $M$, real analytic
  in $\Omega$, and let $\tilde S$ be defined similarly. Consider
  \begin{equation}
    V = \{x\in \Omega: \rho(x)>0,\  \tilde{\rho}(x)<0\} \cup \{x\in \Omega:
   \rho(x)<0,\ \tilde{\rho}(x)>0\}.
  \end{equation}
  It is enough to prove that either $\overline{V}=\Omega$ or
  $V=\varnothing$. Suppose that $V\ne \varnothing$. If $\overline
  V\ne\Omega$, then there exists $q\in (\partial V\cap \Omega) \setminus
  (S\cup \tilde S)$. Since both $\rho$ and $\tilde{\rho}$ change the
  sign near $q$, we conclude that $q$ is an interior point of
  $\overline{V}$. This shows that $\overline{V}=\Omega$.
\end{proof}

\subsection{Cones: real and hermitian signatures}
\label{realconessubsection}

Let $M$ be a quadratic cone in $\R^N$, i.e., $M=\{\rho=0\}$, where
$\rho:\R^N\rightarrow\R$ is a homogeneous quadratic polynomial, and
$M$ is irreducible of dimension $N-1$. Since the Weierstrass
pseudo-polynomial that represents $M$ has degree two, the projection
from $M$ to any $N-1$ dimensional linear subspace is always a
two-sheeted covering on some open set near the origin. It follows
that a quadratic cone can never be represented as the graph of a
continuous function.

After a linear change of coordinates on $\R^N$, and changing $\rho$
to $-\rho$ if required, we can assume that
\begin{equation}\label{diagonalized}
 \rho(x) =\sum_{0<i\leq p}x_i^2 -\sum_{p<i\leq p+q}x_i^2,
\end{equation}
  where $(p,q)$ will be called  the {\em real signature} of $M$ with
$p\geq q\geq 1$. ( Note that if $q=0$, $M$ reduces to a point, and
if $p=q=1$, $M$ is reducible.) The real geometry of cones with $q=1$
and $q>1$ show certain differences: if $q>1$, the set $\R^N\setminus
M$ has two components, and $M^{\rm reg}$ is connected. If $q=1$,
then $\R^N\setminus M$  has {\em three} components. If $\rho$ is
given by \eqref{diagonalized}, and $\Omega^\pm$ are as in
Lemma~\ref{omegas}, the set $\Omega^-$ consists of two components
whereas the set $\Omega^+$ is connected. Also, $M^{\rm reg}$ is
disconnected and has two components.

We now consider real quadratic cones in complex space. Let $\rho(z)$
be  a real-valued homogeneous polynomial of degree two in $\C^n$,
and let $M$ be a quadratic cone given by \eqref{cone-def}. We
associate with  $M$ a pair of non-negative integers $(\pi,\nu)$,
where $\pi$ (resp. $\nu$) is the number of positive (resp. negative)
eigenvalues of the hermitian form
\[
 h_\rho(z,w) = \sum_{i,j=1}^n
\frac{\partial^2\rho}{\partial\overline{ z_i}
\partial {z_j}}\overline{z_i} {w_j} = \tps{\overline{z}} H {w},
\]
where $\tps{A}$ denotes the transpose of a matrix $A$. We call
$(\pi,\nu)$ the {\em hermitian signature} of $M$.  Since $\rho=0$
and $-\rho=0$ define the same cone $M$, we can always assume that
$\pi\geq \nu$.

Let $t_\rho(z)= \rho(z)-h_\rho(z,z)$. Then for each $i,j$,
\[\frac{\partial^2t_\rho}{\partial \overline{z_i}
\partial {z_j}}\equiv 0,\] i.e., $t_\rho$ is a real pluriharmonic
homogeneous polynomial of degree two in $\C^n$. There exists a
holomorphic homogeneous  polynomial $s_\rho$ of degree two in
$\C^n$, such that $t_\rho(z)=\Re(s_\rho)$, and therefore,
\begin{equation}\label{decomposition}
\rho(z) =\Re(\tps{z} S z) + \tps{\overline{z}} H z,
\end{equation}
where $S\in \Sym(n,\C)$, the space of symmetric complex matrices,
and $H$ is an $n\times n$ hermitian matrix. We refer to $\Re(\tps{z}
S z)$ and $\tps{\overline{z}} H z$ as the {\em harmonic} and the
{\em hermitian} parts of the form $\rho$ respectively.

\subsection{Tangent cones}
For a set $E\subset \R^n$, with $0\in E$, the {\em tangent cone}
$T_0 E$ to $E$ at the origin is defined as the set of all limit
 vectors $v\in\R^n$ for all sequences of vectors of the form $t_j
a_j$, where $a_j\in E$, $\lim_{j\to\infty}a_j= 0$, and $t_j>0$.
Clearly, if $E$ is a smooth manifold near the origin, then $T_0 E$
coincides with the tangent plane to $E$ at the origin (hence the
notation). In general, $T_0 E$ is a real cone in the sense that
$tv\in T_0 E$ for $v\in T_0 E$, and $t>0$. We note that if $M$ is a
real quadratic cone in $\R^n$ (or $\C^n$), then $T_0 M=M$. Indeed,
if $a\in M$, then $ta\in M$, and therefore, $ta\in T_0 M$ for all
$t\in\R$. In particular, $a\in T_0 M$. On the other hand, if $a\in
T_0 M$, then there exists a sequence of point $a_j\in M$, $a_j\to
0$, and $t_j>0$ such that $a=\lim_{j\to\infty} t_j a_j =a$. Since
$t_ja_j\in M$, and $M$ is closed, we conclude that $a\in M$.

\begin{lemma}
\label{cone1} Let $M,M'\subset\C^n$ be two germs of quadratic cones
at 0 which are biholomorphic. Then we can take the biholomorphism to
be a complex linear map.
\end{lemma}

\begin{proof}
Let $F$ be a biholomorphic map from $M$ to $M'$. Then by
\cite[Prop.1, Sec. 8.2]{chirka}, $dF_0(T_0M)=T_0 M'$, where $dF_0$
is the differential of $F$ at 0. Since $M=T_0M$ and $M'=T_0 M'$, we
see that $M,M'$ are in fact linearly biholomorphic.
\end{proof}

It follows that the real and hermitian signatures are not only
linear but in fact {\em biholomorphic} invariants of germs of
quadratic cones at the origin.

\section{CR manifolds and functions.}
\label{crdefsection}
\subsection{Definitions and examples}
Recall that an (embedded) CR manifold is a ($\mathcal{C}^1$-) smooth
manifold in $\C^n$ such that the dimension of $H_p M = T_p M \cap i
T_p M$ is independent of $p\in M$. The (complex) dimension of $H_p
M$ is then called the CR dimension of $M$. In particular, any smooth
real hypersurface in $\C^n$ is a CR manifold of CR dimension $n-1$.

Let $M$ be a   smooth real hypersurface in $\C^n$.  Recall that a
distribution $f$ on $M$ is said to be {\em CR} if it satisfies the
tangential Cauchy-Riemann equations. If further $f\in L^1_{\rm
loc}(M)$, and $M$ is orientable, this means that for any
$C^\infty$-smooth $(n,n-2)$-form $\phi$ with compact support on $M$
one has
\begin{equation}\label{cr:definition}
  \int_{M} f \overline{\partial} \phi =0.
\end{equation}

Suppose now that $M$ is an irreducible real analytic hypersurface in
some domain $\Omega$ in $\C^n$. Near every  $p\in M^{{\rm reg}}$,
$M$ is a CR manifold of CR dimension $n-1$. Furthermore (see e.g
\cite{federer}), $M$ defines a  degree one current  on $\Omega$,
denoted by $[M]$, which acts on a compactly supported test form
$\phi$ of degree $2n-1$ on $\Omega$ by
\begin{equation}
  [M](\phi)=\int_{M^{{\rm reg}}}\phi.
\end{equation}
We write  $[M]=[M]^{0,1}+ [M]^{1,0}$ for  the natural splitting of
$[M]$ into a sum of currents of bidegree (0,1) and (1,0). Let
$d\lambda_{2n-1}$ denote the $(2n-1)$-dimensional Lebesgue measure
on $M^{\rm reg}$. We say that $f\in L^1_{\rm loc}(M)$, if for any
compact set $K\subset\Omega$, we have $\int_{M^{{\rm reg}}\cap
K}|f|d\lambda_{2n-1}<\infty$. Then $f[M], f[M]^{0,1}$ and
$f[M]^{1,0}$ are well-defined currents. If $M^{{\rm reg}}=M$, then
\eqref{cr:definition} simply means that $f$ is CR whenever the
current $f[M]^{0,1}$ is $\overline\partial$-closed in $\Omega$. This
leads to the following

\begin{definition}\label{cr:def}
  Let $M$ be a real analytic hypersurface (possibly with
  singularities) in a domain $\Omega\subset\C^n$,
  and $f\in L^1_{\rm loc}(M)$. Then $f$ is called CR if
  for any $C^\infty$-smooth $(n,n-2)$-form $\phi$ with compact support
  on $\Omega$,
  $$\int_{M^{\rm reg}} f \overline\partial \phi=0,$$
  or, equivalently, if $\overline\partial (f[M]^{0,1})=0$.
\end{definition}

It is important to note that according to this definition, a
function $f$ in  $L^1_{\rm loc}(M)$, which is CR on $M^{\rm reg}$,
is not necessarily CR on $M$, even when $M^{\rm sng}$ is a single
point, see e.g., \cite[Example 2.2]{kytmanov}. In other words, the
singularity of $M$ is not in general CR-removable for CR functions
on $M^{\rm reg}$. However, if we assume further conditions regarding
the growth of $f$ near $M^{\rm sng}$, we may conclude that $f$ is CR
on $M$. The following result is also proved in \cite{kytmanov}.  Let
$M=\{\rho(z)=0\}$ be a real analytic hypersurface, and suppose that
$\sigma=M\cap \{d\rho=0\}$ has $2n-1$ measure zero. Let $\tilde
\rho$ be a smooth non-negative function vanishing precisely on
$\sigma$, and let $M_\epsilon = \{z\in M: \tilde \rho> \epsilon\}$.
If $f\in L^1_{\rm loc}(M)$ is CR on $M\setminus \sigma$ and
satisfies the condition
\begin{equation}\label{cr:sufficient}
  \int_{\partial M_\epsilon \cap K} |f(z)| d\lambda_{2n-2} = o(1),\
 {\rm as\ }\epsilon\to 0,
\end{equation}
for all compacts $K$ in $\Omega$, then $f$ is a CR function (as
stated in Definition \ref{cr:def}). In fact this result can be
extended to the case when $\rho$ is merely smooth. This may be used
to construct examples of CR functions on quadratic cones.

\begin{lemma}\label{crexamples}
(a) Let $M\subset\C^n$ be a real quadratic cone , and let $f\in
L^1_{\rm loc}(M)$ be CR on $M^{\rm reg}$. If $f$ is bounded in a
neighbourhood of $M^{\rm sng}$, then $f$ is  CR on $M$.

(b)  Suppose further that $M^{\rm sng}=\{0\}$ and
\[ f(z)=O\left(\frac{1}{\abs{z}^\alpha}\right),\]
where $\alpha<2n-2$ . Then $f$ is   CR on $M$.
\end{lemma}

\begin{proof}
(a) $M^{{\rm sing}}=\{\rho=0, d\rho=0\}$ is a real linear subspace
of $\C^n$ of real dimension at most $2n-2$. If $\dim M^{{\rm
sing}}<2n-2$ and $f$ is bounded, then \eqref{cr:sufficient} is
clearly satisfied, and $f$ is CR on $M$. Suppose now that $\dim
M^{{\rm sing}}=2n-2$. After a suitable $\R$-linear change of
variables in $\C^n=\R^{2n}$ we may assume that $\{d\rho=0\} =
\{x_1=x_2=0\}$. It is easy to see that after an additional change of
variables, the defining function of $M$ admits the form $\rho =
x_1^2 - Ax_2^2$, $A\in \R$. Thus $M$ reduces to the union of real
hyperplanes, but this is ruled out by our assumption. This proves
the assertion.

(b) We can take $\tilde{\rho}(z)=\abs{z}^2$. Since
$\lambda_{2n-2}(\partial M_\epsilon)=O(\epsilon^{2n-2})$, equation
 \eqref{cr:sufficient} holds.
\end{proof}

\subsection{The jump formula and proof of Corollary~\ref{cr-ext}}
\label{jumpsubsection} The {\em jump formula} (\cite{ah},
\cite{chirka:cr})  for CR functions generalizes the classical
Sohocki\v{i}-Plemelj formula for functions of one variable.  If
$\Omega$ is a domain in $\C^n$, with $H^{0,1}(\Omega)=0$ (e.g.,
$\Omega$ can be taken to be pseudoconvex), $M$ is a  smooth
hypersurface in $\Omega$ which divides $\Omega$ into two connected
components $\Omega^+$ and $\Omega^-$, and $f$ is a CR function on
$M$, then there exist functions $F^+\in \mathcal{O}(\Omega^+)$ and
$F^-\in \mathcal{O}(\Omega^-)$ such that the following holds:
\begin{equation}\label{jump}
f=F^+ - F^-.
\end{equation}
The latter equality is understood in an appropriate sense depending
on the smoothness of $f$. In particular, if $f\in L^1_{\rm loc}$,
then for any point $p\in M$, there exists a neighbourhood $U$ such
that
\begin{equation}
  \lim_{\epsilon\to 0^+} \int_{M\cap U} |F^+(\zeta +\epsilon\nu(\zeta)) -
  F^-(\zeta - \epsilon\nu(\zeta)) -f(\zeta)| d\lambda_{2n-1}(\zeta) =0,
\end{equation}
 and $\nu(\zeta)$ is the unit normal vector to~$M$ at
$\zeta\in M$. And if $f$ is a H\"{o}lder continuous function  on $M$
of class $\mathcal{C}^\alpha$, $0<\alpha<1$, then the functions
$F^+$ and $F^-$ extend to $M$ as $\mathcal{C}^\alpha$ functions.

Tr\'epreau's theorem  \cite{trepreau:theorem} now states that in the
situation above, if $M$ is a $C^2$-smooth minimal hypersurface, then
for any point $z\in M$ there exists a neighbourhood $U$ such that
for any CR function $f$ on $M$ either $F^+$ or $F^-$ in \eqref{jump}
extends holomorphically to $U$. Thus any CR function admits a {\em
one-sided} holomorphic extension.

We now consider an analog of the jump formula for real analytic
hypersurfaces. If $\Omega$ is a domain in $\C^n$, with
$H^{0,1}(\Omega)=0$, $M$ is a real analytic hypersurface in
$\Omega$, and $f$ a CR function on $M$ then by Dolbeault's theorem
the $\overline\partial$-closed current $f[M]^{0,1}$ of bidegree
(0,1) is $\overline\partial$-exact. Hence, there exists a current of
degree zero, i.e., a distribution $F$ such that $\overline\partial F
=f[M]^{0,1}$. Now $M$ divides $\Omega$ into two parts $\Omega^\pm$
as in Lemma~\ref{omegas} above. (We emphasize again that $\Omega^+$
and $\Omega^-$ need not be connected.) Since $f[M]^{0,1}$ has
support on $M$,  we have $\overline\partial F =0$ on
$\Omega\setminus M$, and it follows that $F^\pm:=F|_{\Omega^\pm} \in
\mathcal O(\Omega^\pm)$. As in \cite{chirka:cr}, a
Bochner-Martinelli type integral formula can be used to deduce the
correspondence between the values of $F^\pm$ and $f$ at smooth
points of $M$. This yields the jump formula \eqref{jump} for CR
functions defined on (singular) real analytic hypersurfaces. In this
case, \eqref{jump} holds only on $M^{\rm reg}$, in a sense depending
on the smoothness of $f$.

As an immediate  consequence of the jump formula, we can deduce
Corollary~\ref{cr-ext} from Theorem~\ref{extensionthm}. Let $M$ be a
quadratic cone in $\C^n$ which does not have two-sided support at
the origin, and let $f$ be a CR function in a neighbourhood of 0 in
$M$. By the jump formula, we write $f=F^+-F^-$, with
$F^\pm\in\mathcal{O}(\Omega^\pm)$. By Theorem~\ref{extensionthm},
one of the functions $F^\pm$ extends to a neighbourhood $U$ of 0 in
$\C^n$. For definiteness assume that $F^+$ extends to
$\tilde{F}^+\in\mathcal{O}(U)$. Then we can take
$F=\tilde{F}^+-F^-$. For the last statement of
Corollary~\ref{cr-ext}, note that by Lemma~\ref{crexamples}, the
function $f$ is CR on $M$, and therefore, as above, $F=\tilde
F^+-F^-$ is the desired one-sided extension. The equality
$F|_{M^{\rm reg}}=f$ follows from \cite{rosay}. This completes the
proof of the corollary.

Combining the jump formula with the results in \cite{chirka:graphs}
one may further conclude that if $M$ is a real analytic hypersurface
that can be represented as the graph of a continuous function over a
suitable domain in $\C^{n-1}\times \R$, then all CR functions on $M$
extend to one side of $M$, unless $M$ is non-minimal. However, as we
saw in Section~\ref{realconessubsection}, a quadratic cone can never
be represented as a graph of a continuous function.  Therefore,
Chirka's result cannot be used. In fact, the presence of a single
non-smooth point on $M$ leads to new phenomena in the behaviour of
CR functions, as illustrated in the introduction. Thus a different
approach is needed to prove  Theorem~\ref{extensionthm}.

\section{Facts from Linear Algebra.}

The following results in Linear Algebra will be used in the
subsequent sections. We denote by $\Sym(2,\R)$ the space of $2\times
2$ real symmetric matrices, by $SL(2,\R)$ (resp. $SL(2,\C)$) the
group of $2\times 2$  real (resp. complex) matrices of  determinant
1, and by $SO(1,1)$ the group of $2\times 2$ real matrices $g$ with
 $\det g=1$, such that $g^t A g =A$, where $A=\dg(1,-1)$ is the diagonal matrix
with the diagonal entries $1$ and $-1$.

\begin{lemma}\label{lemmaA}
  \label{A} (a) A matrix $k$ preserves the hermitian form
  $\Im(z_1\overline z_2)$ if and only if there exists $g\in
  SL(2,\R)$ and $\theta\in\R$ such that $k=e^{i\theta}g$.

  (b) If $P\in\Sym(2,\R)$,
   $P\ne 0$, then there exists
  $g\in SL(2,\R)$ such that $\tps{g}Pg$ has one of the following forms:

  (i) $\pm\sqrt{\det P}\:I_2$, if $\det P>0$ (here $I_2$ is the
  $2\times 2$ identity matrix);

  (ii) $\pm\sqrt{-\det P}\:\dg(1,-1)$, if $\det P<0$;

  (iii) $\dg(1,0)$, if $\det P=0$.
\end{lemma}

\begin{proof} (a) Denote
  \[ E=\left(
  \begin{array}{cc}
    0 & 1 \\
    -1 & 0 \\
  \end{array}
  \right).
  \]
  Then $\Im(z_1\overline{z_2})= \frac{1}{2i}\tps{\overline{z}} E z$.
  Therefore, we have, $\tps{\overline{k}}\left(\frac{1}{2i} E\right) k
  =\frac{1}{2i} E$, so that taking determinants, $\abs{\det k}^2=1$.
  Consequently, we can write $k=e^{i\theta}g$, where $g\in SL(2,\C)$,
  $\theta\in\R$. Then, $\tps{\overline{g}} E g =E$, or
  $\tps{\overline{g}}E=Eg^{-1}$. From this it follows that $g\in
  SL(2,\R)$. This completes the proof of part (a).

  (b) There exists $g\in GL(2,\R)$ such that
\[ \tps{g} P g = \begin{cases}
\pm I_2, &\text{ if $\det P>0$}\\
\pm \dg(1,-1), &\text{ if $\det P<0$}
\end{cases}
\]

  If $\det g>0$, we can write $g=\delta h$, where $\delta\in\R$, and
  $h\in SL(2,\R)$. Then $\tps{h} P h=\pm\frac{ 1}{\delta^2}I_2$.  If
  $\det g<0$, we can write $g=\delta J h$, where $\delta>0$,
  \begin{equation}\label{j}
    J=\left(
    \begin{array}{cc}
      0 & 1 \\
      1 & 0 \\
    \end{array}
    \right),
  \end{equation}
  and $h\in SL(2,\R)$. In either case, we have
  \[ \tps{h}Ph =\begin{cases}
    \pm\sqrt{\det P}I_2, &\text{ if $\det P>0$}  \\
    \pm\sqrt{-\det P} \dg(1,-1), &\text{if $\det P<0.$}
  \end{cases}
\]

  If $\det P=0$ but $P\not =0$, then there is a $g\in GL(2,\R)$ such that
  $\tps{g}Pg= \dg(1,0)$. We can write $g={\rm
  diag}(1,\det(g))h$, if $\det g>0$, and $g=\dg(1,\det(g))Jh$
  if $\det g<0$, where $h\in SL(2,\R)$, and $J$ is as in
  \eqref{j}. Then $\tps{h} P h=\dg(1,0)$ or ${\rm
  diag}(0,1)$.
\end{proof}

\begin{lemma}\label{B}  Let $Q\in\Sym(2,\R)$, $\det Q\leq 0$ and the sum of the four entries
of $Q$ is different from 0. Then
  there is $k\in SO(1,1)$ such that $\tps{k} Q k$ has at least one
  zero entry on the main diagonal.
\end{lemma}

\begin{proof} Let $\R^+$ be the group of positive reals under
multiplication. Let $\sigma(\tau)=\tau +\tau^{-1}$,
$\delta(\tau)=\tau-\tau^{-1}$, and define
\[ \phi(\tau)= \frac{1}{2}\left(
           \begin{array}{cc}
             \sigma(\tau) & \delta(\tau) \\
             \delta(\tau) & \sigma(\tau) \\
           \end{array}
         \right).\]
Then $\phi:\R^+\rightarrow SO(1,1)$ is a group isomorphism. Assume
that
\[ Q=\left(
                 \begin{array}{cc}
                   p & q \\
                   q& r \\
                 \end{array}
               \right),
\]
and let $Q'= \tps{\phi(\tau)} Q \phi(\tau)$. Denoting the
corresponding elements of $Q'$ by primed letters, we compute
\begin{equation}\label{sigmadelta}
\begin{cases}
            p' & =  \frac{1}{4}(\sigma^2p+2\sigma\delta q +\delta^2r) \\
            q' & =  \frac{1}{4}(\sigma\delta p+(\sigma^2+
            \delta^2)q+\sigma\delta r)\\
            r' & = \frac{1}{4}(\delta^2p+2\sigma\delta q+\sigma^2r).
          \end{cases}
\end{equation}
We need to show that either $p'(\tau)=0$ or $q'(\tau)=0$ has a
solution. After simplifications, we obtain
\[\begin{cases}

(p+2q+r)\tau^2 +2(p-r)\tau^2 +(p-2q+r)&=0\\
 (p+2q+r)\tau^2 -2(p-r)\tau^2 +(p-2q+r)&=0.\end{cases}\]
The possible solutions are contained in the formula:
\[ \tau^2 =\frac{\pm(p-r)\pm2\sqrt{q^2-pr}}{p+2q+r}.\]
By hypothesis the denominator is non-zero, and it is easy to see that at least one of the four values of the right
hand side is a positive real number.
\end{proof}

\section{Normal form for quadratic cones in $\C^2$: Proof of
  Theorem~\ref{normalform}.}

Recall that by Lemma~\ref{cone1}, two quadratic cones which are
biholomorphically equivalent are linearly equivalent, and thus cones
with different hermitian signatures $(\pi,\nu)$, $\pi\geq\nu$, are
not biholomorphic and may be considered separately. The only
possibilities are $(\pi,\nu)=(2,0),(1,1),(1,0)$ and $(0,0).$

  \medskip
  $\boldsymbol{(\pi,\nu)=(2,0)}$. After a
complex linear change of coordinates we may assume that the matrix
$H$ in \eqref{decomposition} is $I_2$. According to  \cite[Theorem
II]{takagi:1924} there is a $2\times 2$ unitary matrix $u$, such
that $u^tSu=\dg(A,B)$, where $A,B$ are the singular values of the
matrix $S$, i.e., the non-negative square roots of the eigenvalues
of $\overline{S}^tS=\overline{S}S$. (This is a special case   of
``Singular Value Decomposition" of matrices.)

 Therefore, there exists a  linear change of variables
such that in the new coordinates the defining equation of $M$ takes
the form
  \begin{equation}
    \rho(z)=\Re(Az_1^2+Bz_2^2)+ |z_1|^2+|z_2|^2.
  \end{equation}
  Permuting $z_1$ and  $z_2$ if required, we have $0\le B\le
  A$. Also, $A>1$, since  $M$ must have dimension~$3$.

If two cones of type $\mathcal{M}_{(2,0)}$ with parameters $(A,B)$
and $(A',B')$ are linearly equivalent, by \cite[Theorem
II]{takagi:1924}  the diagonal matrices $\dg(A,B)$ and $\dg(A',B')$
must have the same singular values. Since $0\leq B\leq A$ and $0\leq
B'\leq A'$, this means that $A=A'$ and $B=B'$. This proves the
uniqueness of the normal form.

  \medskip
  $\boldsymbol{(\pi,\nu)=(1,1)}$. Here, after a linear change of
  variables, we may assume that the defining equation of $M$ has the
  form
  \begin{equation}\label{oneone}
  \rho(z)=\Re(\tps{z} S z) + \Im(z_1\overline z_2),
  \end{equation}
   where $S\in
  \Sym(2,\C)$. If $\det S \ne 0$, we may further make a change of
variables of the form $z^*=e^{i\theta}z$,  so
  that $\det S>0$. (Such a change of variables preserves the hermitian part
 $\Im(z_1\overline{z_2})$.)  Let $S=P+iQ$, where $P=(p_{ij})$ and $Q=(q_{ij})$ are both in
  $\Sym(2,\R)$. Then, a computation shows that $\det S=\det P-\det Q +
  i(q_{11}p_{22}+p_{11}q_{22}-2q_{12}p_{12})$.
  Since $\det S\geq 0$,
  \begin{equation}\label{detess}
    \det S=\det P - \det Q, \end{equation} and
  \begin{equation}\label{noimaginary}
    q_{11}p_{22}+p_{11}q_{22}-2q_{12}p_{12}=0.
  \end{equation}
  We consider three cases depending on the sign of $\det P$.

  {\em Case 1.} Suppose that $\det P>0$. By Lemma~\ref{lemmaA}(b)
  there exists $g\in SL(2,\R)$ such that $\tps{g}P g=\pm\sqrt{\det
  P}I_2$. Since $\tps{g}Qg\in\Sym(2,\R)$, there is $k\in SO(2,\R)\subset SL(2,\R)$ such that
  $\tps{k}(\tps{g}Q g)k=\dg(\lambda_1,\lambda_2)$. Since
  $gk\in SL(2,R)$,  by Lemma~\ref{lemmaA}, the change of variables $z=gkz^*$
 leaves $\Im(z_1\overline{z_2})$ invariant, and in the new coordinates we have
  \begin{equation}\label{newdefining}
  \rho(z^*) =\Re(z^{*t} S' z^*) +\Im(z_1^*\overline z^*_2),
  \end{equation}
  where $S'=\pm
  \sqrt{\det P}I_2 +i\dg (\lambda_1,\lambda_2)$. It follows
  from \eqref{noimaginary} that $\pm \sqrt{\det P}
  (\lambda_1+\lambda_2)=0$. Consequently,
  $\lambda_1=-\lambda_2=\sqrt{-\det Q}$, and therefore, $S'={\rm
  diag}(A,\overline A)$, where $A=\sqrt{\det P}+i\sqrt{-\det
  Q}$. Note that $A\ne 0$, and $\Re A, \Im A\ge 0$.
   Hence, $M$ is equivalent to a cone of type $\mathcal{M}_{(1,1)}^2$.

  {\em Case 2.} Suppose that $\det P<0$. Then by Lemma~\ref{lemmaA}(b)
  there exists $g\in SL(2,\R)$ such that $\tps{g}Pg$ is the matrix of
  the form $\pm\sqrt{-\det P}\:\dg(1,-1)$. Since $\det Q=\det P
  -\det S\le 0$, and thanks to \eqref{detess} and \eqref{noimaginary} above, it follows that the sum of the entries of
  $\tps{g}Pg$ cannot be zero, so by Lemma~\ref{B} there exists a matrix $k\in SO(1,1)$
  such that for $Q'=\tps{k}(\tps{g}Qg)k=(q'_{ij})$, where either
  $q'_{11}$ or $q'_{22}$ vanishes. Since the subgroup $SO(1,1)\subset
  SL(2,\R)$ preserves $P'=\tps{g}Pg$, it follows from
  \eqref{noimaginary} that either $q'_{11}\sqrt{-\det P}=0$, or
  $q'_{22}\sqrt{-\det P}=0$. In either case we deduce that
  $q'_{11}=q'_{22}=0$. After the change of variables $z=kgz^*$ the
  defining equation of $M$ has the form \eqref{newdefining}, with
  \begin{equation}\label{sprime} S' = \left(
          \begin{array}{cc}
           a & i\mu \\
            i\mu & - a \\
          \end{array}
        \right),
  \end{equation}
  where $a=\pm\sqrt{-\det P}$, and $\mu=\pm\sqrt{-\det Q}$. Since
  $\det S= \det P-\det Q \ge 0$, it follows that
  $\abs{\mu}\ge \abs{a}$. After an additional change of variables
  \begin{equation}\label{chofvar}
    \left\{\begin{array}{ccc}
    z^*_1 & = & z_1+iz_2 \\
    z^*_2 & = & -(iz_1+z_2)
    \end{array}
    \right.
  \end{equation}
  in the $(z_1,z_2)$ coordinates, the equation of $M$ becomes:
  \begin{equation}\label{m111}
    \begin{array}{c}
      \Re[a(z_1+iz_2)^2 - 2i\mu (z_1+iz_2)(iz_1+z_2) -a(iz_1+z_2)^2]\\
      +\Im((z_1+iz_2)(i\overline z_1-\overline z_2)) \\
      = 2\Re[(a+\mu)z_1^2 + (\mu-a)z_2^2] + \abs{z_1}^2-\abs{z_2}^2.
    \end{array}
\end{equation}
  We make a further linear change of variables $z_1=\eps_1 z_1^*,
  z_2 =\eps_2 z_2^*$, where the $\eps_j$ are $i$ or 1. It is clear that
  the $\eps_j$ can be so chosen that $M$ is given by
  \begin{equation}\label{bigAB} \Re(Az_1^2+Bz_2^2)
  +\abs{z_1}^2-\abs{z_2}^2=0\end{equation} with
  $A=2\abs{a+\mu},B=2\abs{\mu-a}$. Interchanging $z_1$ and $z_2$ if
  required, we get $\mathcal{M}_{(1,1)}^1$.

{\em Case 3.} Now let $\det P=0$, and let the sum of entries of $\tps{g}Pg$ be nonzero. First suppose that $\det Q<0$.
Then there exists $g\in SL(2,\R)$ such that
$\tps{g}Qg=\pm\sqrt{-\det Q}{\rm
  diag}(1,-1)$. By Lemma~\ref{B} there exists $k\in SO(1,1)$ such that
$\tps{k}(\tps{g}Pg)k$ becomes
\[ \left(
     \begin{array}{cc}
       \alpha' & 0 \\
       0 & 0 \\
     \end{array}
   \right)~\mbox{~or~}\left(
                        \begin{array}{cc}
                          0 & 0 \\
                          0 & \gamma' \\
                        \end{array}
                      \right).
\]
Note that the non-diagonal entries vanish since $\det P=0$. From
this and \eqref{noimaginary} we deduce that $P=0$. It follows then
that $M$ is equivalent to a cone of the form (with $\mu\in\R$)
\[ \Re (i\mu(z_1^2-z_2^2))+\Im(z_1\overline{z_2})=0.\]
We again make the change of variables ~\eqref{chofvar}. In the
$(z_1,z_2)$ coordinates, the cone is given by
\[ 2\Re[i\mu(z_1^2-z_2^2)] + \abs{z_1}^2 -\abs{z_2}^2=0,\]
which after another linear change of coordinates reduces to form
$\mathcal{M}_{(1,1)}^1$ with $A=B=\abs{\mu}$.

Now suppose that $\det Q=0$. If $P\ne 0$, then there exists a $g\in
SL(2,\R)$ such that $\tps{g}Pg=\dg(1,0)$. Let
$Q'=\tps{g}Qg=(q'_{ij})$. Then from \eqref{noimaginary} we conclude
that $q'_{22}=0$. Since $\det Q'=0$, $q'_{12}=0$, and $Q'=\dg
(q'_{11},0)$. Thus we may assume that the defining equation of $M$
has the form
\begin{equation}\label{badcase}
  \rho(z)=\Re(az_1^2)+ \Im z_1\overline z_2 =0,
\end{equation}
where $a=1+iq'_{11}$. Let $\alpha$ be a square root of $a$. After a
change of variables $z_1=\alpha z_1^*, z_2=\frac{1}{\alpha}z_2^*$,
we arrive at $\mathcal{M}_{(1,1)}^3$. If $P=0$, $Q\ne 0$, a similar
argument shows that $M$ can be given by \eqref{badcase} with $a\ne
0$, and again we get $\mathcal{M}_{(1,1)}^3$. Finally, $P=Q=0$
corresponds to $\mathcal{M}_{(1,1)}^1$ with $A=B=0$.

 {\em Case 4.} After reducing the matrix $Q$ to the diagonal form
$q {\rm diag}(1,-1)$, $q\in\mathbb R\setminus\{0\}$, assume that in
the new coordinates the entries of $P$ satisfy
$p_{11}+p_{22}+2p_{12}=0$. It then follows from (17) that
\[ P= p\left(
     \begin{array}{cc}
       1 & -1 \\
       -1 & 1 \\
     \end{array}
   \right),
\]
where $p\in \mathbb R$. Thus the equation of the cone becomes
$$ \Re((p+iq)z_1^2-2pz_1z_2 + (p-iq)z_2^2)+\Im(z_1\overline{z_2})=0.$$
If $p=0$, this reduces to the form considered in Case 3, i.e.,
$\mathcal{M}_{(1,1)}^1$. If $q<0$, we make the coordinate change
$z^*=iz$. If now $p<0$, we make the additional change of coordinates
$z_1^*=z_2$, $z_2^*=z_1$, so that we have $p>0,\ q>0$. Finally, we
make the coordinate change  $ z_1^*= {\sqrt{p}}(z_1-z_2), z_2^*=
\frac{1}{2\sqrt{p}}(z_1+z_2)$, which reduces the cone to type
$\mathcal{M}^4_{(1,1)}$, with $A=2q$.

We now consider uniqueness of the normal forms of hermitian
signature $(1,1)$.   Let $M_1,M_2$ be biholomorphic cones with
$(\pi,\nu)=(1,1)$, such that each is given in one of the normal
forms of Theorem~\ref{normalform}. If the harmonic part of $M_j$ is
$\Re(\tps{z}S_j z)$, then it is easy to see that rank of the
matrices $S_1$ and $S_2$ are equal, in particular, either (1) for
$j=1,2$, $\det S_j\not=0$  or (2) $\det S_1=0=\det S_2$. We consider
case (1) first.
\begin{lemma}\label{C} Suppose that the same quadratic cone $M$ is represented in
    two linear coordinate systems by:
    \begin{equation}\label{s1}
     \Re(\tps{z} S_1 z) + \Im(z_1\overline{z_2})=0,
     \end{equation}
    and
    \begin{equation}\label{s2}
    \Re(\tps{z} S_2 z) + \Im(z_1\overline{z_2})=0,
    \end{equation}
    where each $S_j\in\Sym(2,\C)$ is such that  $\det
    S_j>0$. Then
\begin{romanlist}
\item $\det S_1=\det S_2$, and
\item If $S_j= P_j+i Q_j$, $P_j,Q_j\in \Sym(2,\R)$, then
        $\det P_1=\det P_2$ and $\det Q_1=\det Q_2$.
\end{romanlist}
\end{lemma}

\begin{proof} Let $k\in GL(2,\C)$ be the matrix defining the
change of variables between the coordinates in which equations
\eqref{s1} and \eqref{s2} represent the same cone $M$. Since the
hermitian part is invariant under a linear change of coordinates,
 by Lemma~\ref{lemmaA}(a)  $k=e^{i\theta}g$, where $g\in SL(2,\R)$,
 and $\theta\in \R$. Therefore, $S_2 = e^{2i\theta}(\tps{g} S_1 g)$,
so that $\det S_2= e^{4i\theta} \det S_1$. By the positivity of
$\det S_j$ we have  $e^{4i\theta}=1$  and $\det S_2= \det S_1$.

We see that $\epsilon=e^{i\theta}$ is a fourth root of unity.
Therefore,
\begin{align*}
P_2+i Q_2 &=\tps{k} (P_1+iQ_1) k\\
& = \epsilon^2 \tps{g}(P_1+iQ_1)g\\
&= \pm (\tps{g} P_1 g + i\tps{g} Q_1 g).
\end{align*}
Since the entries of $g\in SL(2,\R)$ are real, it follows that  $P_2
=\pm \tps{g} P_1 g$ and $Q_2 =\pm \tps{g} Q_1 g$, and on taking
determinants, we have $\det P_2 = \det P_1$ and $\det Q_2=\det Q_1$.
\end{proof}

First suppose that both $M_1$ and $M_2$ are of type
$\mathcal{M}_{(1,1)}^2$, and are given by \eqref{s1} and \eqref{s2}
respectively, where $S_j=\dg(A_j,\overline{A_j})$ with $\Re A_j>0$
and $\Im A_j\geq 0$. By Lemma~\ref{C}, $\det P_1= (\Re A_1)^2=(\Re
A_2)^2=\det P_2$, and $\det Q_1=-(\Im A_1)^2=-(\Im A_2)^2 =\det
Q_2$. Therefore, $A_1=A_2$, and $M_1$ and $M_2$ have the same normal
form.

Now suppose that $M_1$ is of type $\mathcal{M}_{(1,1)}^1$ and given
by $\Re (a z_1^2 + b z_2^2) + |z_1|^2 -  |z_2|^2=0$ with $0< a\le
b$. We make the change of variables
\begin{align*}
     z_1^* & =  \frac{1}{2}(z_1+iz_2) \\
     z_2^* & =  \frac{1}{2i}(z_1-iz_2),
   \end{align*}
so that the equation of $M_1$ takes the form \eqref{oneone} with
\[ S= \left(
        \begin{array}{cc}
          a-b & i(a+b) \\
          i(a+b) & b-a \\
        \end{array}
      \right).
\]
Therefore, $\det S=4ab>0$. Writing $S=P+iQ$, we have $\det P=
-(b-a)^2$ and $\det Q= - (a+b)^2$. If $M_2$ is of type
$\mathcal{M}_{(1,1)}^2$, then $S_2 =P_2+i Q_2= {\rm
  diag}(A,\overline{A})$ with $\det P_2 =(\Re A)^2$ and $\det Q_2
=-(\Im A)^2$. Since $\Re A>0$, it follows from Lemma~\ref{C} that
$M_2$ is not biholomorphic to $M_1$. Therefore, $M_2$ is also of
type $\mathcal{M}_{(1,1)}^1$. Assume that $M_2$ is given by
$\Re(a'z_1^2+b'z_2^2)+\abs{z_1}^2-\abs{z_2}^2=0$, with $0\le a'\le
b'$. Then $-(b'-a')^2=-(b-a)^2$ and $-(a'+b')^2=-(a+b)^2$, hence,
$a'=a$ and $b'=b$.

For case (2) suppose that $M_j$, $j=1,2$ are cones of hermitian
signature (1,1) in normal form with harmonic parts
$\Re(\tps{z}S_jz)$, where $\det S_1=0=\det S_2$. By inspection of
the table of normal forms, $M_j$ are of type
$\mathcal{M}_{(1,1)}^3$, or of type $\mathcal{M}_{(1,1)}^1$, with
$B=0$. In the latter case we can write the equation of the cone as
\begin{equation}\label{a}
 \Re(a \tps{z}Ez) + \Im(z_1\overline{z_2})=0,
\end{equation}
where $a=\frac{1}{2}A$, and
\[ E=\left(
       \begin{array}{cc}
         1 & i \\
         i & -1 \\
       \end{array}
     \right).
\]
Note that no cone of type $\mathcal{M}_{(1,1)}^1$ is linearly
equivalent to the cone $\mathcal{M}_{(1,1)}^3$. Indeed, if not, then
there would exist $g\in SL(2,\R)$, $\theta\in\R$ and $a\geq 0$ such
that $e^{2i\theta}\tps{g}\dg(1,0)g=a E$. Let
$e^{2i\theta}\tps{g}\dg(1,0)g=(p_{jk})$, then
$\frac{p_{11}}{p_{22}}$ is a positive real number. However, on the
right hand side, the ratio of the diagonal entries is -1, which is a
contradiction.

Now we show that if $a\not=a'$, the cones \eqref{a} and
\begin{equation}\label{aprime}
 \Re(a' \tps{z}Ez) + \Im(z_1\overline{z_2})=0\end{equation}
are not linearly equivalent. If they were equivalent, there would
exist $g\in SL(2,\R)$ and $\theta\in\R$ such that
$e^{2i\theta}\tps{g}(aE)g=a'E$. Assuming without loss of generality
that $a\not =0$, we have (with $\mu=\frac{a'}{a}e^{-2i\theta}$)
\[\tps{g} E =\mu E g^{-1}.\] From this we easily conclude that $g$ is
the identity, and therefore, $\mu=1$, and $a=a'$.

The uniqueness of cones of type $\mathcal{M}^4_{(1,1)}$ follows the same way
as  for types $\mathcal{M}^1_{(1,1)}$ and $\mathcal{M}^2_{(1,1)}$, using Lemma~\ref{C}.
The uniqueness of
normal forms with hermitian signature (1,1) follows.

For later use we note the following.
\begin{corollary}\label{cor1}
For the cone $M=\{\rho=0\}$, where $\rho$ is as in  \eqref{oneone},
assume that $\det S\geq\frac{1}{4}$, and $\det P<0$. Then $M$ is of
type $\mathcal{M}_{(1,1)}^1$ with $A\ne B$ and  $A\geq1$.
\end{corollary}

\begin{proof} It was shown in {\it Case 2} above that such cones can
  be reduced to type $\mathcal{M}_{(1,1)}^1$ with $A\ne B$.
  Using the same notation, we have
  \begin{align*}AB&=4\abs{\mu^2-a^2}\\
    &=4\abs{\det S'}, \mbox{where $S'$ is as in
      equation~\eqref{sprime}}\\
    &=4\abs{\det S}, \mbox{by Lemma~\ref{C}(i)}\\
    &\geq 1.
  \end{align*}
  Thus at least one of $A$ and $B$ must be greater than or equal to
  1. Since $A\geq B$, it follows that $A\geq 1$.
\end{proof}

\medskip
$\boldsymbol{(\pi,\nu)=(1,0)}$. Such a cone is given by $\rho=0$,
where $\rho(z_1,z_2)= \Re[Az_1^2 +2Bz_1z_2 + C z_2^2] +
\abs{z_1}^2$. If $C\not=0$, then
\[ Az_1^2+ 2B z_1 z_2 +  Cz_2^2 = \left(A-\frac{B^2}{C}\right)z_1^2
+ C\left( z_2+\frac{B}{C}z_1\right)^2.\] Therefore, if we set
$w_1=e^{i\theta_1}z_1$ and $w_2= e^{i\theta_2}\sqrt{C} \left(z_2
+\frac{B}{C}z_1\right)$, we obtain form $\mathcal{M}_{(1,0)}^1$ for
an appropriate choice of $\theta_1$ and $\theta_2$. If $C=0$ and
$B\not=0$, we make the change of variables
\begin{align*}
  z_1^* & =  z_1 \\
  z_2^* & =  2Bz_2+ Az_1,
\end{align*}
thus reducing the cone to type $\mathcal{M}_{(1,0)}^2$. If $C=B=0$,
the defining function is of the form $\Re(Az_1^2)+\abs{z_1}^2$. But
then $\dim M < 3$ if $\abs{A}\leq 1$,  and $M$ is reducible if
$\abs{A}>1$, which violates the assumption.

We now show that the normal forms are unique. If $M$ is non-minimal,
it cannot be of type $\mathcal{M}_{(1,0)}^1$. Suppose that the cone
$\Re(Az_1^2+z_2^2)+\abs{z_1}^2=0$ is biholomorphic to the cone
$\Re(A'z_1^2+z_2^2)+\abs{z_1}^2=0$. Then the linear biholomorphism
between the two cones must map the hermitian form $\abs{z_1}^2$ to
itself, and thus must be of the form
\begin{align*}
  z_1^* & =  e^{i\theta}z_1 \\
  z_2^* & =  pz_1+qz_2.
\end{align*}
It follows that $A'=A$.

\medskip
$\boldsymbol{(\pi,\nu)=(0,0)}$. Writing the defining function as
$\Re(q)=0$, where $q$ is a complex quadratic form, we can
diagonalize $q$ by a change of coordinates to obtain  $q=z_1^2$ or
$q=z_1^2+z_2^2$, depending on the rank of $q$. However,
$\Re(z_1^2)=x_1^2-y_1^2$ is reducible, and we are left with type
$\mathcal{M}_{(0,0)}^1$.

\section{Extension  from quadratic cones: proof of
  Theorem~\ref{extensionthm}.}

If (b) holds, then $M$ does not have two-sided support at 0. Indeed,
if not, then there exist complex analytic hypersurfaces
$A^{\pm}=\{\phi^\pm=0\}\subset \overline{\Omega^{\pm}}$. After
shrinking $\Omega$ if necessary, the function $1/\phi^{\pm}$ is
holomorphic in $\Omega^{\mp}$ and does not extend to the origin.
This proves (b)$\Rightarrow$(a).

The proof of the other direction consists of two steps. First we
will consider the case $n=2$. In the next step, we will consider
cones in $\C^n$, $n>2$, and reduce it to the $n=2$ case by a slicing
argument. In both cases the extension will be derived from the
Kontinuit\"{a}tssatz. We note here that in all cases the constructed
family of analytic discs is not contained in $\Omega^-$ when
$\Omega^-$ consists of two components.

\subsection{(a)\boldmath{$\Rightarrow$(b)} for
\boldmath{$n=2$}}\label{n2}
 After a linear change of coordinates, we
may assume that $M=\{\rho=0\}$ is in the normal form of
Theorem~\ref{normalform}. After rescaling the coordinates, we may
further assume that $\Omega=\{z\in\C^n\colon \abs{z}<2\}$. Let
$\B^2$ denote the unit ball in $\C^2$.

$\boldsymbol{{\mathcal M}_{(2,0)}}$. Suppose
$M=\{\Re(Az_1^2+Bz_2^2)+|z_1|^2+|z_2|^2=0\}$. Consider
\begin{equation}\label{lepsilon}
D_\epsilon=\{z\in\B^2\colon A z_1^2 + B z_2^2=\epsilon\}.
\end{equation}
Let $(z_1,z_2)\in D_\epsilon$, with $\epsilon\geq 0$, then
\begin{equation*}
\rho(z_1,z_2)=\epsilon+\frac{1}{A}\abs{\epsilon-Bz_2^2}+\abs{z_2}^2
\ge 0,
\end{equation*}
with equality iff $z_1=z_2=\epsilon=0$. Thus, $D_0\cap M=\{0\}$ and
for $\epsilon>0$, we have $D_\epsilon\subset\Omega^+$. Thus any
$F^+\in\mathcal{O}(\Omega^+)$ extends to a neighbourhood of $0$.

$\boldsymbol{{\mathcal M}_{(1,1)}^1}$. We consider several cases:
\begin{enumerate}

\item $B\leq A\leq 1$. In this case, $M$ has two-sided support at 0,
 since $\rho|_{\{z_2=0\}}\geq 0$ and $\rho|_{\{z_1=0\}}\leq 0$.

\item $B< 1 \le A$. In this case we let
$D_\epsilon=\{(z_1,z_2)\in\B^2\colon z_1=i\epsilon\}$. Then $M\cap
D_0=\{0\}$, and if $\epsilon>0$, we have
\begin{equation*}
\rho(i\epsilon,z_2) = \Re(Bz_2^2) +\epsilon^2(1-A) -\abs{z_2}^2 <0,
\end{equation*}
which shows that $D_\epsilon\subset \Omega^-$.

\item $1\leq B\leq A.$ If $A=B$, then $M$ is clearly
non-minimal. Therefore we assume that $B<A$.
\[ D_\epsilon=\{ z\in\B^2\colon A z_1^2 + B z_2^2=\epsilon\}.\]
Let $\epsilon\geq 0$, and $(z_1,z_2)\in D_\epsilon$. Then,
$z_1^2=\frac{1}{A}(\epsilon - Bz_2^2)$, so that
\[ \rho(z_1,z_2)= \epsilon
+\frac{1}{A}\abs{\epsilon - B z_2^2}-\abs{z_1}^2.\] Now,
\begin{align*}
\frac{1}{A^2}\abs{\epsilon-Bz_2^2}^2-(\epsilon-\abs{z_2}^2)^2&
=-\left(1-\frac{1}{A^2}\right)\epsilon^2 -\left(1-\frac{B^2}{A^2}\right)\abs{z_2}^4\\
&-2\epsilon\left[\left(1-\frac{B}{A^2}\right)x_2^2+
\left(1+\frac{B}{A^2}\right)y_2^2\right] \\
&\leq 0 \end{align*}
 with strict inequality unless
$\epsilon=z_1=z_2=0$. So $D_0\cap M=\{0\}$, and
$D_\epsilon\subset\Omega^-$ if $\epsilon>0$.
\end{enumerate}

For later use we note the following consequence of the above
arguments and Corollary~\ref{cor1}.
\begin{corollary}\label{cor2}
Let $M$ be a cone defined by $\{\rho=0\}$, where $\rho$ is as in
equation~\eqref{oneone}, and $\det S \geq\frac{1}{4}$. Let $S=P+iQ$,
with $P,Q$ real. If $\det P<0$, then either every function in
$\mathcal{O}(\Omega^+)$ or in $\mathcal{O}(\Omega^-)$ extends to a
neighbourhood of the origin.
\end{corollary}

$\boldsymbol{{\mathcal M}_{(1,1)}^2}$. Let $\lambda_1$ be a real
number such that \[ e^{2i\lambda_1}=-\frac{A}{\overline{A}},\] and
let $\lambda_2=\lambda_1+\pi$. Note that $\lambda_j$ are not integer
multiples of $\pi$. Consider the complex lines
$\Lambda_j=\{z\in\C^2\colon z_2=e^{i\lambda_j}z_1\}$ for $j=1,2$.
Then
\begin{equation*}
\rho(z_1,e^{i\lambda_j}z_1)=
\Re[Az_1^2+\overline{A}e^{2i\lambda_j}z_1^2]+\Im(z_1\cdot
e^{-i\lambda_j}\overline{z_1}) =-\abs{z_1}^2 \sin\lambda_j.
\end{equation*}
Since $\sin\lambda_2=\sin(\lambda_1+\pi)=-\sin\lambda_1$, it follows
that $\Lambda_j$ lie in two different sets of $\Omega^\pm$.
Therefore, $M$ admits two-sided support.

$\boldsymbol{{\mathcal M}_{(1,0)}^1}$. We have
\begin{equation*}
\rho(z)=\Re(A z_1^2 +  z_2^2) + \abs{z_1}^2 = (1+A)x_1^2+(1-A)y_1^2
+x_2^2 - y_2^2.
\end{equation*}
Let $D_\epsilon = \{z\in\B^2\colon A z_1^2 + z_2^2=\epsilon\}.$ If
$\epsilon>0$,  for $z\in D_\epsilon$, we have $\rho(z) = \epsilon
+\abs{z_1}^2\geq 0$, with strict inequality unless
$\epsilon=z_1=z_2=0$. Therefore, $D_0\cap M=\{0\}$, and
$D_\epsilon\subset \Omega^+$, if $\epsilon>0$.

$\boldsymbol{{\mathcal M}_{(1,0)}^2, \mathcal{M}^1_{(0,0)}}$. For
these types, $M$ is non-minimal.

\subsection{(a)\boldmath{$\Rightarrow$(b)} for  $n>2$.}\label{ngeq3}
We first note that if $M\subset\C^n$, $n\geq 3$, is a quadratic
cone, and $L$ is a two dimensional complex linear subspace of $\C^n$
such that $M'=M\cap L$ (viewed as a cone in $L$) does not have
two-sided support, then all functions from one side of $M$ extend to
a neighbourhood of the origin. Therefore, for $M$, which does not
have two-sided support, it suffices to construct a two dimensional
subspace $L\subset\C^n$ so that $M\cap L$ does not have two-sided
support. We assume that $M$ is given by \eqref{decomposition} and
consider different hermitian signatures $(\pi,\nu),\ \pi\geq \nu$.

\medskip

\noindent $\boldsymbol{\pi\ge 2}$. Without loss of generality,
\begin{equation}\label{herm}
  \tps{\overline {z}} H z = |z_1|^2 + |z_2|^2 + \sum_{j=3}^n \eps_j |z_j|^2,
  \ \ \eps_j\in\{-1,0,1\},
\end{equation}
and
\begin{equation}\label{harm}
  \tps{z} S z = \Re(Az_1^2+Bz_2^2+az_2z_3+bz_3^2+cz_1z_3 +s(z)),
\end{equation}
where $A,B\ge 0$, $a,b,c\in\C$, and $s(z)$ does not contain
monomials depending on $z_1, z_2$ and $z_3$ only. We consider three
possibilities below.

(i) Either $A>1$ or $B>1$. Let
\begin{equation}\label{slice}
L=\{z\in \C^n:z_j=0, \ j=3,4,\dots,n\}
\end{equation}
be the two-dimensional linear subspace of $\C^n$, and let $M'=M\cap
L$. Then as in the case ${\mathcal M}_{(2,0)}$ above, we obtain a
sequence of analytic discs $\{D_\epsilon\}_{\epsilon>0}$ in
$L\setminus M'\subset\C^n\setminus M$ with the limit disc $D_0$
touching $M'$ (and therefore $M$) at the origin.

(ii) If $A,B \le 1$, and $AB<1$, we consider $L$ as in
\eqref{slice}. Then the set $M'=M\cap L$ has dimension at most one.
Clearly, there exists a sequence of analytic discs
$\{D_{\epsilon}\}_{\epsilon\ge 0}\subset L$, such that
$D_{\epsilon}\cap M'=\varnothing$, for $\epsilon>0$, and $D_{0}\cap
M'=\{0\}$.

(iii) Suppose now $A=B=1$, then $\dim M'=2$. In this case we
slightly change the subspace $L$. Note that the defining function of
$M$ cannot depend on $z_1$ and $z_2$ only, since otherwise, $\dim M
< 2n-1$. Assume first that not all $\eps_j$ in \eqref{herm} equal
zero, say, $\eps_3\ne 0$. If in \eqref{harm}, $a\ne 0$, we consider
\begin{equation}\label{slice1}
L=\left\{z\in \C^n:z_4 =\dots = z_n=0,\ z_3=\alpha z_2\right\},
\end{equation}
with $\alpha>0$. Then
\begin{equation}\label{ML1}
M\cap L = \left\{\Re\left(z_1^2+(1+\alpha a + \alpha^2b)z_2^2 +
cz_1z_2\right) + |z_1|^2 + (1+\eps_3\alpha^2)|z_2|^2=0\right\}.
\end{equation}
We claim that for any choice of $a,b$ and $c$, there exist
$\alpha>0$, arbitrarily close to zero, and a linear change of
variables in $(z_1,z_2)$ in which the equation of $M\cap L$ becomes
\begin{equation}
  \Re(A'z_1^2+B'z_2^2)+ |z_1|^2+|z_2|^2 = 0,
\end{equation}
and either $A'$ or $B'$ is different from 1 (then our result will
follow from parts (i) and (ii) above.) Indeed, since the hermitian
part  is positive definite, the harmonic part is always
diagonalizable. Therefore, we only need to verify that $A'$ and $B'$
are not both 1. After a linear change of variables, \eqref{ML1} can
be written in the form
\begin{equation}
  \Re(\tps{z}S'z) + |z_1|^2 + |z_2|^2 = 0,
\end{equation}
where
\[ S'=\left(
     \begin{array}{cc}
       1 & \frac{1/2c\alpha}{\sqrt{1+\eps_3\alpha^2}} \\
       \frac{1/2c\alpha}{\sqrt{1+\eps_3\alpha^2}} & \frac{1+a\alpha+b\alpha^2}{1+\eps_3\alpha^2} \\
     \end{array}
   \right).\]
There exists a unitary transformation $u$ such that $\tps{u}S'u$ has
the diagonal form $\dg(A',B')$. To prove that $\alpha$ may be chosen
such that at least one of the coefficients on the diagonal is not 1,
it is enough to show that the determinant of $S'$ is not equal to 1
in absolute value. We have
\begin{equation}
  \det S' = \frac{1+a\alpha+(b-1/4c^2)\alpha^2}{1+\eps_3\alpha^2}.
\end{equation}
Since $a\ne 0$, clearly, $|\det S'|$ cannot equal 1 for all $\alpha
>0$ as $\alpha \to 0$. Thus we are in the situation of case (i) or
(ii), and the extension to the origin follows.

Similarly, if in \eqref{harm} $c\ne 0$, we consider the slice
\begin{equation}\label{slice2}
L=\left\{z\in \C^n:z_4 =\dots = z_n=0,\ z_3=\alpha z_1,\
\alpha>0\right\}.
\end{equation}
and $M\cap L$, in the diagonal form, will not have both coefficients
of the harmonic part equal 1 for a suitable choice of $\alpha$.

Suppose now that $a=c=0$ and $b$ is arbitrary. Then consider the
slice $L$ as in \eqref{slice1} but with $\alpha\in
\C\setminus\{0\}$. We have
\begin{equation}
M\cap L = \Re\left(z_1^2+(1+b\alpha^2)z_2^2 \right)
  +|z_1|^2+(1+\eps_3|\alpha|^2)|z_2|^2 .
\end{equation}
In this case, the determinant of the corresponding matrix $S'$
equals $\frac{1+b\alpha^2}{1+\eps_3|\alpha|^2}$, which is not equal
to 1 in absolute value for a suitable choice of $\alpha$.

The remaining case is when $\eps_j=0$ for $j=3,\dots,n$. Then the
exists $z_k$, $k\ne 1,2$, such that the harmonic part of $M$
contains a monomial involving $z_k$. Without loss of generality we
may assume $k=3$, and $M$ is defined by the equation
\begin{equation}
\Re\left(z_1^2+z_2^2+az_2z_3+bz_3^2+cz_1z_3 +s(z)\right) + |z_1|^2 +
|z_2|^2 = 0,
\end{equation}
where at least one of $a,b$ and $c$ is non-zero. A similar analysis
as above shows that for any choice of $a$, $b$ and $c$, there exists
a slice $L$ of the form \eqref{slice1} or \eqref{slice2} such that
$M\cap L$ is a quadratic cone in $\C^2$ with positive definite
hermitian part such that the at least one of the coefficients of the
harmonic part (when reduced to the diagonal form) is different from
1. Thus we obtain the holomorphic extension to the origin. This
proves the extension when the hermitian part of $M$ has at least two
positive eigenvalues.

\medskip

\noindent $\boldsymbol{\pi=\nu=1}.$ Denote $z\in\C^n$ as
$z=(\tilde{z},z')$, where $\tilde{z}=(z_1,z_2)\in \C^2$,
$z'=(z_3,\ldots,z_n)\in\C^{n-2}$. Any cone $M$ in $\C^n$ with
$(\pi,\nu)=(1,1)$ can be written after a change of coordinates as
\begin{equation}\label{oneonenge3}
 \rho= \Re\left(
\tps{\tilde{z}}S\tilde{z} + 2 z_1l_1(z')+ 2z_2l_2(z')+ q(z')\right)
+ \Im(z_1\overline{z_2})=0, \end{equation} where $S\in\Sym(2,\C)$,
$l_1,l_2$ are $\C$-linear forms on $\C^{n-2}$ and $q$ is a complex
quadratic form on $\C^{n-2}$. We consider the cases $q\equiv 0$ and
$q\not\equiv 0$ separately.

$\boldsymbol{q\equiv 0}$. We begin by noting the following fact.

\begin{lemma}\label{ge3stdformslemma}
Suppose that $M$ is given by \eqref{oneonenge3} with $q\equiv 0$.
Then there exists a coordinate system in which  $M$ is given by
\begin{equation}\label{ge3stdforms}
\Re(\tps{\tilde{z}}S\tilde{z}+2R(z))+\Im(z_1\overline{z_2})=0,
\end{equation}
where $S\in\Sym(2,\C)$, and $R$ is one of the following:
\begin{romanlist}
\item[{\rm(i)}] $0$,
\item[{\rm(ii)}] $z_1z_3$,
\item[{\rm(iii)}] $z_2z_3$,
\item[{\rm(iv)}] $cz_1z_3+ z_2z_3$, $c\in\C\setminus\{0\}$,
\item[{\rm(v)}] $z_1z_3+z_2z_4$.
\end{romanlist}
\end{lemma}

\begin{proof}
If $l_1=l_2=0$, we get (i). If $l_2=0$ but $l_1\ne 0$, then we can
change coordinates in the $z'$ variable to ensure that $l_1=z_3$,
which gives us (ii). Similarly, if $l_1=0$ and $l_2\ne 0$ then we
get (iii). If $l_1$ and $l_2$ are linearly dependent, but neither is
0, we get (iv).  Finally, if $l_1$ and $l_2$ are linearly
independent, we can change coordinates in the last $n-2$ variables
$z'$ such that $l_1=z_3$, $l_2=z_4$, which gives us (v).
\end{proof}

We will now consider the various possibilities for $R$.

{(i)} $R=0$. Then $M$ is given by
$\Re(\tps{\tilde{z}}S\tilde{z})+\Im(z_1\overline{z_2})=0$. This
defining function involves only the variables $z_1$ and $z_2$. Let
$L$ be the two dimensional  subspace of $\C^n$ given by
$\{z_k=0,k\geq 3\}$. Let $M'=M\cap L$. $M'$ can be identified with
the cone in $\C^2$ given by
$\Re(\tps{z}Sz)+\Im(z_1\overline{z_2})=0$. Clearly $M'$ has
two-sided support, iff $M$ does.  The result follows from the
discussion at the beginning of this subsection
(subsection~\ref{ngeq3}).

{(ii)} $R=z_1z_3$. We let
\begin{equation}\label{abbc}
S= \left(
                            \begin{array}{cc}
                              A & B \\
                              B & C \\
                            \end{array}
                          \right).
\end{equation}
If $C=0$,  the cone $M$ contains the complex hypersurface
$\{z_1=0\}$, and therefore is non-minimal. If $C\ne 0$, we let $L$
be the two dimensional linear subspace of $\C^n$ given by
\begin{equation}\label{alphabetaeqn}
\left\{\begin{array}{cccc}
           z_3 & =&\alpha z_1+\beta z_2 &\\
           z_k & =&0,& k\geq 4.
         \end{array}
\right.
\end{equation}
 Consider the two dimensional slice $M\cap L$, which is given by
  \eqref{alphabetaeqn} and
\[\Re(\tps{\tilde{z}}S_* \tilde{z}) +\Im(z_1\overline{z_2})=0,
\]
where
\[ S_*=\left(
                            \begin{array}{cc}
                              A+2\alpha & B+\beta \\
                              B +\beta& C \\
                            \end{array}
                          \right).\]
If $\Re(C)=0$, we let\[ \begin{array}{ccc}
         \alpha & = & -\frac{1}{2}( A+{C}) \\
         \beta & = & -B+\sqrt{\abs{C}^2+2}
       \end{array}.
\]
Then $\det S_*= -2$. After a change of variables
$z^*=e^{i\frac{\pi}{4}}z$, the equation of the cone in the new
coordinates becomes $\Re(\tps{\tilde{z}}\tilde{S} \tilde{z})
+\Im(z_1\overline{z_2})=0,$ where  $\tilde{S}=iS_*$. Let
$\tilde{S}=P+iQ$, with $P,Q$ real matrices. Then $\det \tilde{S}=2$,
and  $P=\dg(-iC,iC)$, so that $\det P = -(iC)^2<0$. Therefore, we
can apply Corollary~\ref{cor2}. If $\Re(C)\ne 0$, we let
\[ \begin{array}{ccc}
         \alpha & = & -\frac{1}{2}( A+\overline{C}) \\
         \beta & = & -B+i\sqrt{\abs{C}^2+2}
       \end{array}.
\]
Then $\det S_*= 2>1$. If $S_*=P+iQ$, where $P,Q\in \Sym(2,\R)$, then
$P=\linebreak[4] \dg(-\Re C,\Re C)$,  and $\det P=-(\Re C)^2<0$.
Again, we can apply Corollary~\ref{cor2}.

{(iii)} $R=z_2z_3$. This is similar to case (ii).

{(iv)} $R=cz_1z_3+z_2z_3, \ c\in\C\setminus\{0\}$. If
$c\in\R\setminus\{0\}$, the linear change of variables
$(w_1,w_2,w')=(cz_1+z_2,\frac{1}{c}z_2,z')$ reduces the situation to
case {(ii)}. Observe that the matrix corresponding to this change of
coordinates in the $z_1,z_2$ variables belongs to $SL(2,\R)$, and
therefore preserves the hermitian part $\Im(z_1\overline{z_2})$.
(Lemma~\ref{A}(a)).

If $c\in\C\setminus\R$, then $M$ is given by
\begin{equation}\label{c}
\Re(\tps{\tilde{z}}S\tilde{z}+
2(cz_1+z_2)z_3)+\Im(z_1\overline{z_2})=0,
\end{equation}
where $S$ is as in \eqref{abbc}. For $\alpha\in\C\setminus\{0\}$,
let $L_\alpha\subset\C^n$ be the two dimensional subspace given by
\begin{equation}\label{lalpha2} \left\{\begin{array}{ccc}
            z_3 & = & \alpha z_2 \\
            z_k & = & 0, \ k\geq 4.
          \end{array}\right.
\end{equation}
Then $M\cap L_\alpha$ is given by the above equations and the
equation $\Re(\tps{\tilde{z}} (S+\alpha
T)\tilde{z})+\Im(z_1\overline{z_2})=0$, where
\[ T=\left(
             \begin{array}{cc}
               0 & c \\
               c & 2 \\
             \end{array}
           \right).\]
We set $ S_* =S+ \alpha T$ , and let $\zeta$ be such that
\begin{equation}\label{zeta}
\zeta^4 =\frac{\overline{\det S_*}}{\abs{\det S_*}}.
\end{equation}
After a change of coordinates  $z^*=\zeta z$, the cone is
represented by $\Re(\tps{\tilde{z}} \tilde{S}
\tilde{z})+\Im(z_1\overline{z_2})=0$, where $\tilde{S}=\zeta^2 S_*$.
By the choice of $\zeta$, $\det \tilde{S}>0$. We now obtain an
asymptotic expression for $\zeta^2$ as
$\abs{\alpha}\rightarrow\infty$:
\begin{align}
\pm\zeta^2 &= \sqrt{\frac{\overline{A}\overline{C} +
2\overline{A}\overline{\alpha}-
(\overline{B}+\overline{c}\overline{\alpha})^2}
{\abs{AC + 2A\alpha- (B+c\alpha)^2}}}\nonumber\\
&=i\frac{\overline{\alpha}}{\abs{\alpha}}\cdot
\frac{\overline{c}}{\abs{c}}+O\left(\frac{1}{\abs{\alpha}}\right).
\label{asymp}
\end{align}
Therefore, as $\abs{\alpha}\rightarrow\infty$,
\[\pm\alpha \zeta^2=  i
\abs{\alpha}\frac{\overline{c}}{\abs{c}}+O(1),\] and
\begin{equation}\label{est1}
\pm\Re(\alpha\zeta^2)=  \abs{\alpha}\sin(\arg c) +O(1).
\end{equation}
Let $\theta$ be such that $\sin(\arg c)$ and $\sin(\theta+\arg c
-\arg A)$ have different signs. (Note that since
$c\in\C\setminus\R$, $\sin(\arg c) \ne 0$.) Let
$\alpha=\abs{\alpha}e^{i\theta}$, where $\abs{\alpha}$ is large and
will be chosen later. Using \eqref{asymp} again,
\begin{equation}\label{est2}
\pm\Re(A\zeta^2)= \abs{A}\sin(\theta+\arg c - \arg A)
+O\left(\frac{1}{\abs{\alpha}}\right).
\end{equation}

We now show that $\abs{\alpha}$ can be chosen such that
$\det\tilde{S}>1$ and if $\tilde{S}=P+iQ$ (with $P$,$Q$ real), then
$\det P<0$. A computation shows that \begin{equation}\label{p} \det
P=\Re(A\zeta^2)\Re(C\zeta^2)+2\Re(A\zeta^2)\Re(\alpha\zeta^2)
-\left(\Re((B+C\alpha)\zeta^2)\right)^2.
\end{equation}
From \eqref{est1} and \eqref{est2} above, we see that there is $K>0$
such that for $\abs{\alpha}$ large enough
\[\Re(A\zeta^2)\Re(\alpha\zeta^2)<-K\abs{\alpha}.\]
Since the first term in the right in \eqref{p} is bounded by
$\abs{A}\abs{C}$, it follows that we can choose $\abs{\alpha}$ so
large that $\det P<0$. Also $\det \tilde{S}\rightarrow\infty$ as
$\abs{\alpha}\rightarrow\infty$, it follows that we can also make
$\det \tilde{S}>\frac{1}{4}$. The result now follows from
Corollary~\ref{cor2}.

{(v)}~$R=z_1z_3+z_2z_4$. If $A\ne 0$ or $C\ne 0$, the situation can
be clearly reduced to {(ii)} or {(iii)}. Therefore, the only new
situation is when $A=C=0$, i.e.,
\[ S = \left(
         \begin{array}{cc}
           0 & B \\
           B & 0 \\
         \end{array}
       \right).
\]
Let $L\subset\C^n$ be the two dimensional linear subspace  defined
by
\[ \left\{ \begin{array}{ccc}
             z_3 & = & \frac{1}{2}z_1 + \left(-\frac{B}{2}+i\right)z_2 \\
             z_4 & = & \left(-\frac{B}{2}+i\right)z_1-\frac{1}{2}z_2 \\
             z_k & = & 0, \ k\geq 5.
           \end{array}\right.
\]
Then $M\cap L$ is given by the above equations and
\[ \Re(\tps{\tilde{z}} S_* \tilde{z})+ \Im(z_1z_2)=0,
\]
where \[ S_*= \left(
                       \begin{array}{cc}
                         1 & 2i \\
                         2i & -1 \\
                       \end{array}
                     \right).
\]
It follows that $\det S_*=3$, and if $S_*=P+iQ$, then $P=
\dg(1,-1)$, so that $\det P=-1<0$. Again, extension follows from
Corollary~\ref{cor2}.

$\boldsymbol{q\not\equiv 0.}$ In the $\tilde{z}=(z_1,z_2)$ variable,
we make the change \eqref{chofvar}, and in the
$z'=(z_3,\ldots,z_n)$, we make a linear change of variables so that
$M$ is given by
\[\Re\left(
\tps{\tilde{z}}S\tilde{z} + 2 z_1l_1(z')+ 2z_2l_2(z')+
\sum_{j=3}^kz_j^2\right) + \abs{z_1}^2 -\abs{z_2}^2,\] where $l_1$
and $l_2$ are complex linear functionals on $\C^{n-2}$, and $3\leq
k\leq n$. Now let $L$ be given by
\begin{equation}\label{l}z_k=0 ~~~\mbox{~~~for}~~~  k\geq 4.\end{equation}
Then $M\cap L$ is given by the equations \eqref{l} and
\[ \Re(\tps{\tilde{z}}S\tilde{z} + 2c_1 z_1z_3 + 2 c_2 z_2z_3 + z_3^2) +
\abs{z_1}^2 -\abs{z_2}^2=0,
\]
where $c_1,c_2\in \C$. Define\[\left\{
\begin{array}{ccc}
  z_3^* & = & c_1z_1+c_2z_2+z_3\\
  z_k^* & = & z_k, \mbox{~~~~for~~} k\not= 4\\
\end{array}.\right.
\]
With respect to the new coordinates (after suppressing the asterisk
for convenience), $M\cap L$ is given by equations~\eqref{l} along
with
\[ \Re(\tps{\tilde{z}}S'\tilde{z} +  z_3^2) + \abs{z_1}^2 -\abs{z_2}^2=0,
\]
where $S'\in \Sym(2,\C)$. Assume that
\[ S= \left(
        \begin{array}{cc}
          A & B \\
          B & C \\
        \end{array}
      \right),\]
and let $L_\alpha\subset L$ be the  two dimensional linear subspace
of $\C^3$ defined by \begin{equation}\label{lalpha} z_2=\alpha z_1,
\end{equation}
where $\alpha$ is a complex number such that $\abs{\alpha}\ne 1$.
Then $M\cap L_\alpha$ is given by \eqref{l},\eqref{lalpha} and
\[ \Re((A+2B\alpha +C\alpha^2)z_1^2 +  z_3^2) +
(1-\abs{\alpha}^2)\abs{z_1}^2=0.
\]
After rescaling the $z_1$ coordinate, $M\cap L_\alpha$ is
represented by \eqref{l}, \eqref{lalpha}, and
\[\Re\left(\abs{E_\alpha} z_1^2 +  z_3^2\right) +
\abs{z_1}^2=0, \ {\rm where} \ E_\alpha =\frac{A+2B\alpha
  +C\alpha^2}{1-\abs{\alpha}^2}.
\]
Repeating the construction of the family of analytic discs for type
$\mathcal{M}_{(1,0)}^1$, we prove that for this cone holomorphic
functions on one side have extension to the origin.

\medskip

\noindent $\boldsymbol{\pi=1,\ \nu=0.}$ We can write $M$ as
\[ \Re(Az_1^2+z_1l(z')+q(z'))+\abs{z_1}^2=0,\]
where $z'\in\C^{n-1}$, $l$ and $q$ are complex linear and quadratic
forms on $\C^{n-1}$ respectively. Note that if $q\equiv 0$, we have
$\{z_1=0\} \subset M$, so that $M$ is non-minimal. We will therefore
assume that $q\ne 0$.

If $l\equiv 0$, by a linear change of variables in $z'$, we can
assume that $q(z')=\sum_{j=2}^k z_j^2$, where $2\leq k\leq n$.
Taking a slice by $L=\{z_k=0,k\geq 3\}$ we see that $M\cap L$ is of
type $\mathcal{M}_{(1,0)}^1$, which does not have two-sided support.

If $l\not=0$, we can assume that $l(z')=z_2$ after a change of
variables. If $\frac{\partial q }{\partial z_2}\not=0$, then we can
again take a slice by $L=\{z_k=0,k\geq 3\}$, and again we see that
$M\cap L$ is of type $\mathcal{M}_{(1,0)}^1$.

If $\frac{\partial q }{\partial z_2}=0$, we can change variables in
the $(z_3,\ldots,z_n)$ variables to diagonalize $q$ and write
$q(z')=\sum_{j=3}^k z_j^2$, where $3\leq k\leq n$. We take a slice
by $L=\{ z_k=0, k=2, \mbox{~or~} k\geq4\}$. Then $M\cap L$ is a cone
of type $\mathcal{M}_{(1,0)}$ in the $z_1z_3$ plane. This completes
the proof for $\pi=1,\nu=0$.

\medskip

\noindent $\boldsymbol{\pi=\nu=0}$. These cones are clearly
non-minimal. Theorem~\ref{extensionthm} is proved.

\section{Classification of quadratic cones with two-sided support.}
As a corollary to the proof of Theorem~\ref{extensionthm} we obtain
the following complete classification of real quadratic cones in
$\C^n$ which admit two-sided support by complex hypersurfaces.

\begin{proposition}
Let $M\subset\C^n$ be a real quadratic cone in $\C^n$ which admits
two-sided support at 0. Then one of the following holds.

(i) $M$ is biholomorphic to $M'\times \C^{n-2}$, where $M'$ has the
normal form either of type $\mathcal{M}_{(1,1)}^1$ with $0\le B\le
A\le 1$, or $A=B$, or one of $\mathcal{M}_{(1,1)}^2$,
$\mathcal{M}_{(1,1)}^3$, $\mathcal{M}_{(1,1)}^4$, $\mathcal{M}_{(1,0)}^2$, or
$\mathcal{M}_{(0,0)}^1$.

(ii) $M$ is biholomorphic either to
\begin{equation}\label{ts1}
\{\Re(z_1^2+\cdots+z_k^2)=0,\  2< k\le n \},
\end{equation}
or
\begin{equation}\label{ts2}
\{\Re(z_1z_2+z_1\overline{z_3})=0\}.
\end{equation}
\end{proposition}

\begin{proof}
(i) If $n=2$, this follows immediately from Section~\ref{n2}. If
$n\geq 3$, it is clear that $M=M'\times\C^{n-2}$ has two-sided
support iff $M'$ has two-sided support in $\C^2$, so the result
follows in this case as well.

(ii) An examination of the proof in Section~\ref{ngeq3} shows that
  every cone in three variables, which is not extendable, is, in fact,
  non-minimal. Let $A=\{\alpha=0\}\subset M$ be the germ of an
  irreducible complex hypersurface at the origin. Then, for the
  tangent cones, $T_0A\subset T_0M=M$, but since $A$ is a complex
  hypersurface, $T_0A=\{ \alpha_\mu=0\}$, where $\alpha_\mu$ is the
  homogeneous polynomial consisting of the nonzero terms of the least
  degree in the Taylor expansion of $\alpha$. It follows that we can
  take $\alpha$ to be a homogeneous complex polynomial, and since
  $\rho$ has degree 2, $\alpha$ has degree 1 or 2. Therefore, there is
  a real polynomial $\phi$ on $\C^n$ such that
  $\rho=\phi\alpha+\overline\phi\overline\alpha=\Re(\phi\alpha)$. If
  $\rho$ is of degree 2, $\phi$ has degree 0, and $\rho=\Re(c\alpha)$
  for some  constant $c$. After a linear change of variables, we can
  diagonalize the complex quadratic form $c\alpha$, and obtain
  \eqref{ts1}. If $\deg \alpha=\deg\phi=1$, then
  $\phi=\lambda+\overline{\mu}$, where $\lambda$ and $\mu$ are complex
  linear maps from $\C^n$ to $\C$. First suppose that
  $\{\alpha,\phi,\mu\}$ is  a linearly independent set in the dual a
  space of $\C^n$. We can find a system of coordinates
  $(z_1,\ldots,z_n)$ on $\C^n$ such that $z_1=\alpha$, $z_2=\lambda$
  and $z_3=\mu$. Therefore, we have $\rho=\Re(\phi\alpha) =
  \Re((z_2+\overline{z_3})z_1)$, which gives \eqref{ts2}. If $\alpha$,
  $\lambda$, $\mu$ are linearly dependent, we are reduced to
  classifying non-minimal cones in~$\C^2$.
\end{proof}

\section{Two-sided support for smooth hypersurfaces: proof of
Proposition~\ref{smoothconda}.}

Proposition~\ref{smoothconda} will be deduced from the following
lemma, for which we do not claim any originality. For completeness,
we give a proof essentially following
\cite[Lem.~3.2]{trepreau:theorem}.

\begin{lemma}\label{trepreaulem}
Let $\Omega\subset\C^n$ be a domain whose boundary can be  locally
represented as a Lipschitz graph and let $p\in\partial\Omega$.
Suppose that there is an open set $U\subset\C^n$ containing $p$ and
a holomorphic function $f\in\mathcal{O}(\Omega\cap U)$ which does
not extend to any neighbourhood of $p$, and  let $A$ be a germ at
$p$ of  complex analytic sets such that $A\subset\overline{\Omega}$.
Then $A\subset
\partial\Omega$.
\end{lemma}

\begin{proof} Since the problem is local, we can assume that the
  coordinates $(z',z_n)\in\C^{n-1}\times\C$ have been so chosen that
  $p=0$, and
  \[ \Omega= \{z\in\C^n, (z',x_n)\in U\colon y_n > g(z',x_n)\},\]
  where $g$ is a real-valued Lipschitz function defined in a
  convex neighbourhood $U$ of $0$ in $\C^{n-1}\times\R$.  It
  is known \cite{trepreau:theorem}, \cite{chirka:graphs} that the envelope of
  holomorphy $\widehat{\Omega}$ of $\Omega$ is schlicht  and of the form
  \[ \widehat{\Omega}= \{z\in\C^n,(z',x_n)\in U\colon y_n > \hat{g}(z',x_n)\},\]
  where $\hat{g}\leq g$ is an upper semi-continuous function. Since
  the origin is a point of non-extendability, it follows that
  $\hat{g}(0,0)=g(0,0)=0$.

By a standard abuse of notation, denote by $A$ a particular
representative of the germ $A$ at the origin. Without loss of
generality we may assume that $A$ is irreducible. Arguing by
contradiction, suppose that there exists $w\in A\setminus
\partial\Omega$ arbitrarily close to the origin. Let $A_w$ be an
irreducible complex curve in $A$ connecting $w$ and the origin. Note
that such $A_w$ always exists. Indeed, if $\dim A=1$, then $A_w =
A$, so assume that $\dim A = d \ge 2$. Then near the origin $A$ can
be represented as a branched covering $\pi : A \to L$ over a
$d$-dimensional subspace $L\subset \mathbb C^n$. Further, in some
neighbourhood $V$ of the origin, $\pi^{-1}(0)=0$. We may assume that
$w\in V$. Let $l\subset L$ be the complex line connecting $\pi(w)$
and the origin. Then one of the irreducible components of
$\pi^{-1}(l)$ is the desired curve.

Using a Puiseux parametrization  we obtain a holomorphic map $\phi$
from the closed unit disc $\Delta$ into $\C^n$ such that
$\phi(\Delta)$ is a neighbourhood of $0$ in $A_w$, and
$\phi(\Delta)\subset\overline{\Omega}\cap U\cap V$, where $U$ is a
neighbourhood of 0 as in the statement of the lemma, i.e., some
$f\in\mathcal{O}(\Omega\cap U)$ does not extend to any neighbourhood
of $0$, and $V$ is a neighbourhood of the origin such that there is
a constant $C\geq 1$ such that for $z\in V$,
\begin{equation}
\label{geometric}{\rm dist}(z,\partial\Omega)\geq
C^{-1}(y_n-g(z',x_n)).\end{equation} Such $V$ exists since
$\partial\Omega$ is Lipschitz.

Applying an automorphism of $\Delta$, we can assume that
$\phi(0)=0$. Let $\chi(z)= -\log {\rm  dist}(z,\partial\Omega)$ and
$\hat{\chi}(z)= -\log {\rm dist}(z,\partial\widehat{\Omega})$. Then
$\hat{\chi}$ is a plurisubharmonic exhaustion of $\widehat{\Omega}$,
and on $\Omega$ we have $\chi\leq \hat{\chi}$. For $t> 0$ and $\zeta
\in \Delta$ define
\[ \psi_t(\zeta)= \hat{\chi}(\phi(\zeta)+t ie_n),\]
where $e_n=(0,\dots,0,1)$ is the $n$-th standard basis element of
$\C^n$. Then $\psi_t$ is a subharmonic function for each $t$.
Further, define for $k>0$
\[ \mu_k(t)= {\rm Lebesgue~ measure~of~
}\{\zeta\in\Delta\colon\psi_t(\zeta)<k \}.\] To prove our result it
is sufficient to show that for all $k$ we have  $\lim_{t\rightarrow
0^+}\mu_k(t)=0$. We claim that for small~$t$,
\begin{equation}\label{aa}
\psi_t(0) \geq \log\left(\frac{1}{t}\right),
\end{equation}
and
\begin{equation}\label{ab}
\max_{\zeta\in\Delta}(\psi_t(\zeta))\leq\log\frac{C}{t},
\end{equation}
where $C$ is as in ~\eqref{geometric}.

For \eqref{aa} note that since $0\in\partial\widehat{\Omega}$, we
have
 ${\rm dist}(ite_n,\partial\widehat{\Omega})\leq {\rm dist}(ite_n,0)=t$,
 so $\psi_t(0)=-\log{\rm  dist}(ite_n,\partial\widehat{\Omega})\geq -\log t.$

For \eqref{ab} observe that there is a nonempty subset $T\subset
\partial\Delta$ such that $\phi(T)\subset\partial\Omega$ (otherwise
by the Kontinuit\"{a}tssatz $f$ would extend to a neighbourhood of
0.) From \eqref{geometric},
\begin{align*}
{\rm dist}(z+ite_n,\partial\Omega)&\geq
\frac{1}{C}(y_n+t-g(z',x_n))\\
& =\frac{1}{C}(y-g(z',x_n))+\frac{t}{C} \end{align*} In particular,
if $z\in\partial\Omega$, then ${\rm
dist}(z+ite_n,\partial\Omega)\geq \frac{t}{C}$. Since
$\phi(\zeta)\in \partial\Omega$ for $\zeta\in T$, we have for small
$t$,
\[
\min_{\zeta\in\partial\Delta}{\rm
dist}(\phi(\zeta)+ite_n,\partial\Omega) \geq \frac{t}{C}.\]
Therefore,
\begin{align*}
\max_{\zeta\in\Delta}(\psi_t(\zeta))
&=\max_{\zeta\in\partial\Delta}\hat{\chi}(\phi(\zeta)+ite_n)\\
&\leq\max_{\zeta\in\partial\Delta}\chi(\phi(\zeta)+ite_n)\\
 &= -\log\min_{\zeta\in\partial\Delta}{\rm dist}(\phi(\zeta)+it e_n,\partial\Omega)\\
 &\leq -\log\frac{t}{C},
\end{align*}
which proves the second claim. Now apply the sub-averaging property
to the subharmonic function $\psi_t$:
\begin{align*}
\pi\log\left(\frac{1}{t}\right)&\leq\pi \psi_t(0)\\
&\leq \int_\Delta \psi_t(\zeta)d\xi d\eta\\
&= \int_{\{\psi_t(\zeta)<k\}} + \int_{\{\psi_t(\zeta)\geq k\}}\\
&\leq k\mu_k(t) + (\pi-\mu_k(t) ) \max_{\zeta\in\Delta}\psi_t(\zeta)\\
&\leq  k\mu_k(t) + (\pi-\mu_k(t))\log\left(\frac{C}{t}\right)\\
&=\left(k-\log\left(\frac{C}{t}\right)\right)\mu_k(t) +
\pi\log\left(\frac{C}{t}\right).
\end{align*}
Therefore, for small $t$, we have:
\begin{align*}
\mu_k(t)&\leq\frac{\pi\left[\log\left(\frac{1}{t}\right)-
\log\frac{C}{t}\right]}{k-\log C+\log t}\\
&= \frac{-\pi\log C}{k-\log C+\log t}\\
&\rightarrow 0 {\mbox{~ as ~}} t\rightarrow 0^+.
\end{align*}
\end{proof}
\begin{proof}[Proof of Proposition~\ref{smoothconda}]
Let $A^-$ be defined in a neighbourhood $U$ of $p$ in $\C^n$ by an
equation $\{z\in U:f(z)=0\}$. Then
$\frac{1}{f}\in\mathcal{O}(\Omega^+\cap U)$ does not extend to any
neighbourhood of $p$. Since $A^+\subset\overline{\Omega}$, by
Lemma~\ref{trepreaulem}, $A^+\subset M$. By the same argument,
$A^-\subset M$.  If $A^+\not=A^-$, then $A=A^+\cup A^-$ is a
non-smooth complex hypersurface through $p$, and $A\subset M$. But
any complex hypersurface contained in a hypersurface represented as
a graph is smooth (See  \cite[Lemma~3]{chirka:rado}.)
\end{proof}

\end{document}